\documentclass{amsart}

\newtheorem{theorem}{Theorem}[section]
\newtheorem{lemma}[theorem]{Lemma}

\newtheorem{proposition}[theorem]{Proposition}
\newtheorem{corollary}[theorem]{Corollary}

\theoremstyle{definition}

\theoremstyle{remark}

\usepackage{graphicx}
\usepackage{amsmath}
\usepackage{amsfonts}
\usepackage{amssymb}
\newcommand{\be}{\begin{equation}}

\newcommand{\ee}{\end{equation}}

\newcommand{\Om}{\Omega}

\newcommand{\K}{\Upsilon}

\newcommand{\om}{\omega}

\newcommand{\G}{\Gamma}

\newcommand{\La}{\Lambda}

\newcommand{\si}{\sigma}

\newcommand{\dz}{\wedge}

\newcommand{\ba}{\begin{array}}

\newcommand{\ea}{\end{array}}

\newcommand{\beq}{\begin{eqnarray}}

\newcommand{\eeq}{\end{eqnarray}}

\newtheorem{lm}{lemma}

\newtheorem{thee}{theorem}

\newtheorem{proo}{proposition}

\newtheorem{co}{corollary}

\newtheorem{rem}{remark}

\newtheorem{deff}{definition}

\newcommand{\bd}{\begin{deff}}

\newcommand{\ed}{\end{deff}}

\newcommand{\bl}{\begin{lm}}

\newcommand{\el}{\end{lm}}

\newcommand{\bp}{\begin{proo}}

\newcommand{\ep}{\end{proo}}

\newcommand{\bt}{\begin{thee}}

\newcommand{\et}{\end{thee}}

\newcommand{\bc}{\begin{co}}

\newcommand{\ec}{\end{co}}

\newcommand{\brm}{\begin{rem}}

\newcommand{\erm}{\end{rem}}

\newcommand{\der}{{\rm d}}

\usepackage{t1enc}

\def\Bbb{\mathbb}
\def\Cal{\mathcal}

\newcommand{\newc}{\newcommand}

\let\ccdot\cdot
\def\cdot{\hbox to 2.5pt{\hss$\ccdot$\hss}}

\newc{\aR}{\mbox{\boldmath{$ R$}}}
\newc{\aS}{\mbox{\boldmath{$ S$}}}
\newc{\aT}{\mbox{\boldmath{$ T$}}}
\newc{\aW}{\mbox{\boldmath{$ W$}}}

\newc{\aK}{\mbox{\boldmath{$ K$}}}
\newc{\aL}{\mbox{\boldmath{$ L$}}}

\newcommand{\ce}{{\Cal E}}

\newcommand{\cq}{{\Cal Q}}

\newcommand{\bR}{{\Bbb R}}
\newcommand{\bT}{{\Bbb T}}

\usepackage{amssymb}
\usepackage{amscd}

\newcommand{\nd}{\nabla}

\newcommand{\Rho}{{\mbox{\sf P}}}
\newcommand{\Up}{\Upsilon}

\newcommand{\wh}{\widehat}

\newcommand{\IT}[1]{{\rm(}{\it{#1}}{\rm)}}

\newcommand{\Pa}{{\Bbb I}}

\newcommand{\nn}[1]{(\ref{#1})}

\newcommand{\bg}{\mbox{\boldmath{$ g$}}}

\newcommand{\vol}{\mbox{\large\boldmath $ \epsilon$}}

\let\t=\tau

\let\G=\Gamma

\newcommand{\V}{{\mbox{\sf P}}}                   
\newcommand{\J}{{\mbox{\sf J}}}
\newc{\obstrn}[2]{B^{#1}_{#2}}

\newcommand{\rpl}                         
{\mbox{$
\begin{picture}(12.7,8)(-.5,-1)
\put(0,0.2){$+$}
\put(4.2,2.8){\oval(8,8)[r]}
\end{picture}$}}

\newcommand{\lpl}                         
{\mbox{$
\begin{picture}(12.7,8)(-.5,-1)
\put(2,0.2){$+$}
\put(6.2,2.8){\oval(8,8)[l]}
\end{picture}$}}

\usepackage{ifthen}

\newc{\tensor}[1]{#1}
\newc{\Mvariable}[1]{\mbox{#1}}
\newc{\down}[1]{{}_{#1}}
\newc{\up}[1]{{}^{#1}}

\newc{\JulyStrut}{\rule{0mm}{6mm}}
\newc{\midtenPan}{\mbox{\sf S}}
\newc{\midten}{\mbox{\sf T}}
\newc{\midtenEi}{\mbox{\sf U}}
\newc{\ATen}{\mbox{\sf E}}
\newc{\BTen}{\mbox{\sf F}}
\newc{\CTen}{\mbox{\sf G}}

\def\sideremark#1{\ifvmode\leavevmode\fi\vadjust{\vbox to0pt{\vss
 \hbox to 0pt{\hskip\hsize\hskip1em
 \vbox{\hsize3cm\tiny\raggedright\pretolerance10000
 \noindent #1\hfill}\hss}\vbox to8pt{\vfil}\vss}}}%
                        
\numberwithin{equation}{section}

\begin{document}
\title{Obstructions to conformally Einstein metrics in $n$ dimensions} 
\vskip 1.truecm \author{A.\ Rod Gover} \address{Department of Mathematics,
University of Auckland, Private Bag 92019, Auckland, New Zealand}
\email{r.gover@auckland.ac.nz} \thanks{This research was supported by
the Royal Society of New Zealand (Marsden Grant 02-UOA-108), the
New Zealand Institute for Mathematics and its Applications, and the KBN grant 2 P03B
12724 }

\author{Pawe\l~ Nurowski} \address{Instytut Fizyki Teoretycznej,
Uniwersytet Warszawski, ul. Hoza 69, Warszawa, Poland}
\email{nurowski@fuw.edu.pl} \thanks{During the preparation of this
article P. N. was a member of the VW Junior Research Group ``Special
Geometries in Mathematical Physics'' at Humboldt University in
Berlin. P. N. would also like to thank the University of
Auckland for hospitality during the preparation of
this article.}

\date{April 7, 2004}
\begin{abstract}
We construct polynomial conformal invariants, the vanishing of which
is necessary and sufficient for an $n$-dimensional suitably generic
(pseudo-)Riemannian manifold to be conformal to an Einstein manifold.
We also construct invariants which give necessary and sufficient
conditions for a metric to be conformally related to a metric with
vanishing Cotton tensor. One set of invariants we derive generalises
the set of invariants in dimension four obtained by Kozameh, Newman
and Tod. For the conformally Einstein problem, another set of
invariants we construct gives necessary and sufficient conditions for
a wider class of metrics than covered by the invariants recently
presented by M.\ Listing.  We also show that there is an alternative
characterisation of conformally Einstein metrics based on the tractor
connection associated with the normal conformal Cartan bundle. This
plays a key role in constructing some of the invariants. Also using
this we can interpret the previously known invariants geometrically in
the tractor setting and relate some of them to the curvature of the
Fefferman-Graham ambient metric.
\end{abstract}
\maketitle
\section{Introduction}

The central focus of this article is the problem of finding necessary
and sufficient conditions for a Riemannian or pseudo-Riemannian
manifold, of any signature and dimension $n\geq 3$, to be locally
conformally related to an Einstein metric.  In particular we seek
invariants, polynomial in the Riemannian curvature and its covariant
derivatives, that give a sharp obstruction to conformally Einstein
metrics in the sense that they vanish if and only if the metric
concerned is conformally related to an Einstein metric. For example in
dimension 3 it is well known that this problem is solved by the Cotton
tensor, which is a certain tensor part of the first covariant
derivative of the Ricci tensor. So 3-manifolds are conformally
Einstein if and only if they are conformally flat. The situation is
significantly more complicated in higher dimensions. Our main result
is that we are able to solve this problem in all dimensions and for
metrics of any signature, except that the metrics are required to be
non-degenerate in the sense that they are, what we term, weakly
generic. This means that, viewed as a bundle map $TM\to\otimes^3 TM$,
the Weyl curvature is injective.  The results are most striking for
Riemannian $n$-manifolds where we obtain a single trace-free rank two
tensor-valued conformal invariant that gives a sharp
obstruction. Setting this invariant to zero gives a quasi-linear
equation on the metric.  Returning to the setting of arbitrary
signature, we also show that a manifold is conformally Einstein if and
only if a certain vector bundle, the so-called standard tractor
bundle, admits a parallel section. This powerful characterisation of
conformally Einstein metrics is used to obtain the sharp obstructions
for conformally Einstein metrics in the general weakly generic
pseudo-Riemannian and Riemannian setting.  It also yields a simple
geometric derivation, and unifying framework, for all the main
theorems in the paper.

The study of conditions for a metric to be conformally Einstein has a
long history that dates back to the work of Brinkman \cite{B1,B2} and
Schouten \cite{S}. Substantial progress was made by Szekeres in the
1963 \cite{Sz}. He solved the problem on 4-manifolds, of signature $-2$, by
explicitly describing invariants that provide a sharp
obstruction. However his approach is based on a spinor formalism and
is difficult to analyse when translated into the equivalent tensorial
picture. In the 1980's Kozameh, Newman and Tod (KNT) \cite{CTP} found
a simpler set of conditions.  While their construction was based
on Lorentzian 4-manifolds the invariants obtained provide obstructions
in any signature. However these invariants only give a sharp
obstruction to conformally Einstein metrics if a special class of
metrics is excluded (see also \cite{CTPa} for the reformulation of the
KNT result in terms of the Cartan normal conformal connection). Baston and Mason \cite{Mas}
proposed another pair of conformally invariant obstruction invariants
for 4-manifolds. However these give a sharp obstruction for a smaller
class of metrics than the KNT system (see \cite{BaiEast}).

One of the invariants in the KNT system is the conformally invariant
Bach tensor.  In higher even dimensions there is an interesting higher
order analogue of this trace-free symmetric 2-tensor due to Fefferman
and Graham and this is also an obstruction to conformally Einstein
metrics \cite{FeffGr,GoPetprogress,GrH}. This tensor arises as an
obstruction to their ambient metric construction. It has a close
relationship to some of the constructions in this article, but this is
described in \cite{GoPetprogress}. Here we focus on invariants which
exist in all dimensions.  Recently Listing \cite{Listing} made a
substantial advance. He described a trace-free 2-tensor that gives, in
dimensions $n\geq 4$, a sharp obstruction for conformally Einstein
metrics, subject to the restriction that the metrics are what he terms
``nondegenerate''. This means that the Weyl curvature is maximal rank
as a map $\Lambda^2 TM\to \Lambda^2 TM$. In this paper metrics
satisfying this non-degeneracy condition are instead termed
$\Lambda^2$-generic.

Following some general background, we show in Sections \ref{necessary}
and \ref{KNTstyle} that it is possible to generalise to arbitrary
dimension $n\geq 4$ the development of KNT. This culminates in the
construction of a pair of (pseudo-)Riemannian invariants $F^1_{abc}$
and $F^2_{ab}$ whose vanishing is necessary and sufficient for the
manifold to be conformally Einstein provided we exclude a small class
of metrics (but the class is larger than the class failing to be
$\Lambda^2$-generic). See theorem \ref{II}. These invariants are {\em
natural} in the sense that they are given by a metric partial
contraction polynomial in the Riemannian curvature and its covariant
derivatives. $F^1$ is conformally covariant and $F^2$ is conformally
covariant on metrics for which $F^2$ vanishes. Thus together they form
a conformally covariant system.

In Section \ref{newinvts} we show that very simple ideas reveal new
conformal invariants that are more effective than the system $F^1$ and
$F^2$ in the sense that they give sharp obstructions to conformal
Einstein metrics on a wider class of metrics. Here the broad treatment
is based on the assumption that the metrics are weakly generic as
defined earlier. This is a strictly weaker restriction than requiring
metrics to be $\Lambda^2$-generic; any $\Lambda^2$-generic metric is
weakly generic but in general the converse fails to be true.  One of
the main results of the paper is Theorem \ref{nob} which gives a
natural conformally invariant trace-free 2-tensor which gives a sharp
obstruction for conformally Einstein metrics on weakly generic
Riemannian manifolds. Thus in the Riemannian setting this improves
Listing's results.  In Riemannian dimension 4 there is an even simpler
obstruction, see Theorem \ref{4ob}, but an equivalent result is in
\cite{Listing}. In Theorem \ref{Lam2ob} we also recover Listing's main
results for $\Lambda^2$-generic metrics as special case of the general
setup.  In all cases the invariants give quasi-linear equations.  The
results mentioned are derived from the general result in Proposition
\ref{Eprop}.
 We should point out that while this Proposition does not
in general lead to natural obstructions, in many practical situations,
for example if a metric is given explicitly in terms of a basis field,
this would still provide an effective route to testing whether or not
a metric is conformally Einstein, since a choice of tensor $\tilde{D}$
can easily be described. (See the final remark at the end of Section
\ref{newinvts}.)  

In section \ref{newinvts} we also pause, in Proposition \ref{Cob} and
Theorem \ref{L2Cotton}, to observe some sharp obstructions to metrics
being conformal to a metric with vanishing Cotton tensor. We believe
these should be of independent interest. Since the vanishing of the
Cotton tensor is necessary but not sufficient for a metric to be
Einstein, it seems that the Cotton tensor could play a role in setting
up problems where one seeks metrics suitably ``close'' to being
Einstein or conformally Einstein.

In Section \ref{tractor}, following some background on tractor
calculus, we give the characterisation of conformally Einstein metrics
as exactly those for which the standard tractor bundle admits a
(suitably generic) parallel section. The standard (conformal) tractor
bundle is an associated structure to the normal Cartan conformal
connection.  The derivations in Section \ref{tensor} are quite simple
and use just elementary tensor analysis and Riemannian differential
geometry. However they also appear ad hoc. We show in Section
\ref{tractor} that the constructions and invariants of Section
\ref{tensor} have a natural and unifying geometric interpretation in
the tractor/Cartan framework. This easily adapts to yield new
characterisations of conformally Einstein metrics, see Theorem
\ref{tracchar}. From this we obtain, in Corollary \ref{weakanysig},
obstructions for conformally Einstein metrics that are sharp for
weakly generic metrics of any signature. Thus these also improve on the results in
\cite{Listing}.

We believe the development in Section \ref{tractor} 
should  have an important role in suggesting how an analogous programme
could be carried out for related conformal problems as well as
analogues on, for example, CR structures where the structure and
tractor calculus is very similar. We also use this machinery to show
that the system $F^1$, $F^2$ has a simple interpretation in terms of
the curvature of the Fefferman-Graham ambient metric.

Finally in Section \ref{examples} we discuss explicit metrics to shed
light on the invariants constructed and their applicability. This
includes examples of classes metrics which are weakly generic but not
$\Lambda^2$-generic.
 Also here,
as an example use of the machinery on explicit metrics, we identify the
conformally Einstein metrics among a special class of 
Robinson-Trautman metrics. 

The authors wish to thank Ruibin Zhang, Paul-Andi Nagy and Michael
Eastwood for very helpful discussions.


\noindent

\section{Conformal characterisations via tensors} \label{tensor}

In this section we use standard tensor analysis on (pseudo-)Riemannian manifolds to derive 
sharp obstructions to conformally Einstein metrics.

\rm
\noindent
\subsection{Basic (pseudo-)Riemannian objects} \label{background}

Let $M$ be a smooth manifold, of dimension $n\geq 3$, equipped with a
Riemannian or pseudo-Riemannian metric $g_{ab}$. We employ Penrose's
abstract index notation \cite{ot} and indices should be assumed
abstract unless otherwise indicated. We write $\ce^a$ to denote the
space of smooth sections of the tangent bundle on $M$, and $\ce_a$ for
the space of smooth sections of the cotangent bundle.  (In fact we
will often use the same symbols for the corresponding bundles, and
also in other situations we will often use the same symbol for a given bundle and
its space of smooth sections, since the meaning will be clear by
context.) We write $\ce$ for the space of smooth functions and all
tensors considered will be assumed smooth without further comment.  An
index which appears twice, once raised and once lowered, indicates a
contraction.  The metric $g_{ab}$ and its inverse $g^{ab}$ enable the
identification of $\ce^a$ and $\ce_a$ and we indicate this by raising
and lowering indices in the usual way.

The metric $g_{ab}$ defines the Levi-Civita connection $\nabla_a$ with
the curvature tensor $R^a_{~bcd}$ given by
$$(\nabla_a\nabla_b-\nabla_b\nabla_a)V^c={R_{ab}^{~~~c}}_d V^d 
\quad\text{
where} \quad \ V^c\in \ce^c.$$
This can be decomposed into the totally trace-free {\em Weyl curvature}
$C_{abcd}$ and the symmetric {\em
Schouten tensor} $\Rho_{ab}$ according to
$$
R_{abcd}=C_{abcd}+2g_{c[a}\Rho_{b]d}+2g_{d[b}\Rho_{a]c}.
$$
Thus $\Rho_{ab}$ is a trace modification of the Ricci tensor 
$R_{ab}=R_{ca}{}^c{}_b$:
$$R_{ab}=(n-2)\Rho_{ab}+\J g_{ab},~~~~~~~~~~\J:=\Rho^a_{~a}.$$
Note that the Weyl tensor has the symmetries 
$$
 C_{abcd}=C_{[ab][cd]}=C_{cdab},\quad  C_{[abc]d}=0.
$$
where we have used the square brackets to denote the
antisymmetrisation of the indices. 

We recall that the metric $g_{ab}$ is an Einstein metric if the trace
free part of the Ricci
tensor vanishes. This condition, when written in terms of the
 Schouten tensor, is given by 
$$\Rho_{ab}-\frac{1}{n} \J g_{ab}=0.$$
In the following we will also need
the Cotton tensor $A_{abc}$ and the Bach tensor $B_{ab}$. These are defined by
\begin{equation} 
A_{abc}:=2\nabla_{[b}\Rho_{c]a}\label{cot} 
\end{equation} and 
\begin{equation} B_{ab}:=\nabla^c
A_{acb}+\Rho^{dc}C_{dacb}.\label{bach} 
\end{equation}
It is straightforward to verify that the Bach tensor is symmetric.
From the contracted Bianchi identity $\nd^a\Rho_{ab}=\nd_b \J$ it follows that the 
Cotton tensor is totally trace-free. Using this, and that the Weyl tensor is trace-free, 
it follows that the Bach tensor is also trace-free.

Let us adopt the convention that sequentially labelled indices are
implicitly skewed over. For example with this notation the  
Bianchi symmetry is simply $R_{a_1a_2a_3 b}=0$. Using this symmetry and 
 the definition (\ref{cot}) of $A_{ba_1a_2}$ we obtain a useful identity
\begin{equation}\label{ABianchiid}
\nabla_{a_1}A_{ba_2a_3}=\Rho_{a_1}{}^c C_{a_2a_3bc} .
\end{equation}
Further important identities arise from the Bianchi
identity $\nabla_{a_1}R_{a_2a_3 de}=0$:
\begin{equation} \label{WBianchiid} 
 \nd_{a_1}C_{a_2a_3 cd} =g_{ca_1}A_{da_2a_3}-g_{da_1}A_{ca_2a_3}.
\end{equation} 
\begin{equation} (n-3)A_{abc}=\nabla^d
C_{dabc}\label{bi1} 
\end{equation} 
\begin{equation} 
\nabla^a\Rho_{ab}=\nabla_b \J\label{bi2} 
\end{equation}
\begin{equation} \nabla^a A_{abc}=0.\label{bi3} 
\end{equation} 

\subsection{Conformal properties and naturality}\label{confsect}
Metrics $g_{ab}$ and $\widehat{g}_{ab}$ are said to be conformally
related if 
\begin{equation} \widehat{g}_{ab}={\rm
e}^{2\Upsilon}g_{ab},\label{conftr} \quad \Upsilon\in \ce , 
\end{equation} 
and
the replacement of $g_{ab}$ with $\widehat{g}_{ab}$ is termed a {\em
conformal rescaling}. Conformal rescaling in this way results in a conformal 
transformation of the Levi-Civita connection. This is given by  
\begin{equation}\label{conntrans}
\wh{\nd_a u_b} =\nd_a u_b -\Up_a u_b -\Up_b u_a + g_{ab}\Up^c u_c
\end{equation}
for a 1-form $u_b$. The conformal transformation of the Levi-Civita
connection on other tensors is determined by this, the duality between
1-forms and tangent fields, and the Leibniz rule.

A tensor $T$ (with any number of covariant and
contravariant indices) is said to be {\em conformally covariant} (of
{\em weight} $w$) if, under a conformal rescaling (\ref{conftr}) of
the metric, it transforms according to
$$
T\mapsto\widehat{T}={\rm e}^{w\Up}T,
$$ for some $w\in {\Bbb R}$. We will say $T$ is conformally {\em
invariant}\, if $w=0$. We are particularly interested in natural
tensors with this property. A tensor $T$ is {\em natural} if there is an expression for
$T$ which is a metric partial contraction, polynomial in the metric,
the inverse metric, the Riemannian curvature and its covariant
derivatives.  

The weight of a conformally covariant depends on the placement of
indices.  It is well known that the Cotton tensor in dimension $n=3$
and the Weyl tensor in dimension $n\geq 3$ are conformally invariant
with their natural placement of indices, i.e.\
$\widehat{A}_{abc}=A_{abc}$ and $\widehat{C}_{ab}{}^c{}_d =
{C}_{ab}{}^c{}_d $. In dimension $n\geq 4$, vanishing of the Weyl
tensor is equivalent to the existence of a scale $\Up$ such that the
transformed metric $\widehat{g}_{ab}={\rm e}^{2\Up}g_{ab}$ is flat
(and so if the Weyl tensor vanishes we say the metric is {\em
conformally flat}). In dimension $n=3$ the Weyl tensor vanishes
identically. In this dimension $g_{ab}$ is conformally flat if and
only if the Cotton tensor vanishes. 

An example of tensor which fails to be conformally covariant
is the Schouten tensor.  We have
\begin{equation}\label{Rhotrans}
\Rho_{ab}\to\widehat{\Rho}_{ab}=\Rho_{ab}-\nabla_a\K_b+\K_a\K_b-\frac{1}{2}\K_c\K^c
g_{ab},
\end{equation}
where 
$$
\K_a=\nabla_a\Up.
$$ 
Thus the property of the metric being Einstein is not conformally
invariant. A metric $g_{ab}$ is said to be {\em conformally Einstein}
if there exists a conformal scale $\Up$ such that
$\widehat{g}_{ab}={\rm e}^{2\Up}g_{ab}$ is Einstein.

For natural tensors the property of being conformally covariant or
invariant may depend on dimension. For example it is well known that
the Bach tensor is conformally covariant in dimension 4. In other
dimensions the Bach tensor fails to be conformally covariant.

\subsection{Necessary conditions for conformally Einstein metrics}\label{necessary}

Suppose that $g_{ab}$ is conformally Einstein. As mentioned above this
means that there exists a scale $\Up$ such that the Ricci tensor, or
equivalently the Schouten tensor for $\wh{g}_{ab}:=e^{2\Up} g_{ab}$,
is pure trace. That is
$$
\widehat{\Rho}_{ab}-\frac{1}{n}\widehat{\J}\widehat{g}_{ab}=0.
$$ 
This equation, when written in terms of Levi-Civita connection
$\nabla$ and Schouten tensor $\Rho_{ab}$ associated with $g_{ab}$
reads,
\begin{equation}
\Rho_{ab}-\nabla_a\K_b+\K_a\K_b-\frac{1}{n}Tg_{ab}=0,\label{confein}
\end{equation}
where
$$
T=\J-\nabla^a\K_a+\K^a\K_a.
$$ Conversely if there is a gradient $\Up_a=\nd_a\Up$ satisfying
\nn{confein} then $\wh{g}_{ab}:=e^{2\Up} g_{ab}$ is an Einstein
metric. Thus, with the understanding that $\Up_a=\nd_a\Up$,
(\ref{confein}) will be termed the {\em conformal Einstein equations}.
There exists a smooth function $\Up$ solving these if and only
if the metric $g$ is conformally Einstein.

To find consequences of these equations we apply $\nabla_c$ to both sides of (\ref{confein})
and then antisymmetrise the result over the $\{ca\}$ index pair. Using that
the both the Weyl tensor and the Cotton tensor are completely
trace-free this leads to the first integrability condition which is
$$
A_{abc}+\K^d C_{dabc}=0.
$$
Now taking $\nabla^c$ of this
equation, using the definition of the Bach tensor (\ref{bach}), the
identity (\ref{bi1}), and again this last displayed equation, we get
$$B_{ab}+\Rho^{dc}C_{dabc}-(\nabla^c \K^d -(n-3)\K^d\K^c)C_{dabc}=0.$$
Eliminating $\nabla^c \K^d$ by means of the Einstein condition
(\ref{confein}) yields a second integrability condition:
$$
B_{ab}+(n-4)\K^d\K^c C_{dabc}=0.
$$

Summarising we have the following proposition.
\begin{proposition}
If $g_{ab}$ is a conformally Einstein metric then the corresponding Cotton tensor
$A_{abc}$ and the Bach tensor $B_{ab}$ satisfy the following 
conditions 
\begin{equation} 
A_{abc}+\K^d C_{dabc}=0,\label{1ic} 
\end{equation}
and
\begin{equation}
B_{ab}+(n-4)\K^d\K^c C_{dabc}=0.\label{2ic}
\end{equation}
for some gradient 
$$
\K_d=\nabla_d\Up.
$$ 
\end{proposition}
Here $\Up$ is a function which conformally rescales the metric $g_{ab}$ to an Einstein metric 
$\widehat{g}_{ab}={\rm e}^{2\Up}g_{ab}$. \\

\noindent
{\bf Remarks}
\begin{itemize}
\item
Note that in dimension $n=3$ the first integrability condition
  (\ref{1ic}) reduces to $A_{abc}=0$ and the Weyl curvature
  vanishes. Thus, in dimension $n=3$, if \nn{1ic} holds then \nn{2ic}
  is automatically satisfied and the conformally Einstein metrics are
  exactly the conformally flat metrics. The vanishing of the Cotton
  tensor is the necessary and sufficient condition for a metric to
  satisfy these equivalent conditions.  This well known fact solves
  the problem in dimension $n=3$. Therefore, for the remainder of
  Section \ref{tensor} we will assume that $n\geq 4$.
\item In dimension $n=4$ the second integrability condition reduces to
the conformally invariant Bach equation:
\begin{equation}
B_{ab}=0.\label{bache} 
\end{equation}

\end{itemize}

\subsection{Generalising the KNT characterisation}\label{KNTstyle}

Here we generalise to dimension $n\geq 4$ the characterisation of 
conformally Einstein metrics
given by Kozameh, Newman and Tod \cite{CTP}.
Our considerations are local and so we assume, without loss of
generality, that $M$ is oriented and write $\vol$ for the volume
form.  Given the Weyl tensor $C_{abcd}$ of the metric $g_{ab}$, we
write $C^*_{b_1\cdots b_{n-2}cd}:=\vol_{b_1\cdots
b_{n-2}}{}^{a_1a_2}C_{a_1a_2cd}$. Note that this is completely
trace-free due to the Weyl Bianchi symmetry $C_{a_1a_2a_3 b}=0$.
 Consider the equations
\begin{equation}
C_{abcd}F^{ab}=0\label{c3}, \end{equation} 
\begin{equation} C_{abcd}H^{bd}=0\label{c3bar}, 
\end{equation}
and
\begin{equation}\label{c3*bar}
C^*_{b_1\cdots b_{n-2}cd}H^{b_1d}=0 ,
\end{equation}
for a skew symmetric tensor $F^{ab}$ and a symmetric trace-free tensor
$H^{ab}$. We say that the metric $g_{ab}$ is {\em generic} if and only
if the only solutions to equations (\ref{c3}), (\ref{c3bar}) and
\nn{c3*bar} are $F^{ab}=0$ and $H^{ab}=0$.  Occasionally we will be
interested in the superclass of metrics for which \nn{c3} has only
trivial solutions but for which we make no assumptions about
\nn{c3bar} and \nn{c3*bar}; we will call these {\em
$\Lambda^2$-generic} metrics.  That is, a metric is $\Lambda^2$-generic
if and only if the Weyl curvature is injective (equivalently, maximal
rank) as a bundle map $\Lambda^2 TM\to \Lambda^2 TM$. Let $||C||$
be the natural conformal invariant which is the pointwise determinant
of the map
\begin{equation}\label{Lam2map}
C:\Lambda^2T^*M \to \Lambda^2T^*M,
\end{equation}
given by $W_{ab}\mapsto C_{ab}{}^{cd}W_{cd}$ and write
$\tilde{C}_{abcd}$ for the tensor field which is the pointwise
adjugate (i.e.\ ``matrix of cofactors'') of the Weyl curvature tensor,
viewed as an endomorphism in this way. Then
$$
\tilde{C}_{ef}^{~~ab}C_{ab}^{~~cd}=||C||\delta^{[c}_{~[e}\delta^{d]}_{~f]}
$$
and if $g$ is a $\Lambda^2$-generic metric then $||C||\neq 0$ and we have
\begin{equation}
||C||^{-1}\tilde{C}_{ef}^{~~ab}C_{ab}^{~~cd}=\delta^{[c}_{~[e}\delta^{d]}_{~f]}.\label{til}
\end{equation}
For later use note that it is easily verified that 
$\tilde{C}_{abcd} $ is natural (in fact simply polynomial in
the Weyl curvature) and conformally covariant.

For the remainder of this subsection we consider only
generic metrics, except where otherwise indicated. In this setting, we
will prove that the following two conditions are equivalent:
\begin{itemize}
\item[(i)] The metric $g_{ab}$ is conformally Einstein.
 \item[(ii)] There exists a vector field $K^a$ on $M$ such that the
 following conditions [C] and [B] are satisfied:
$$
~~~~~~~~~A_{abc}+K^d C_{dabc}=0~~~~~~~~~~~~~~~~~~~~~~~~~{\rm [C]}
$$
$$
B_{ab}+(n-4)K^dK^cC_{dabc}=0.~~~~~~~~~~~~~~~~~~~~~{\rm [B]}
$$
\end{itemize}

Adapting a tradition from the General Relativity literature
(originating in \cite{Sz}), we call a manifold for which the metric $g_{ab}$
admits $K^a$ such that condition [C] is satisfied a 
{\em conformal C-space}.

We first note that if a generic metric satisfies condition [C] then
the  field $K_d$ must be a gradient. To see this take $\nabla^a$
of equation [C]. This gives
$$
\nabla^a A_{abc}+C_{dabc}\nabla^aK^d+(n-3)K^aK^dC_{adbc}=0,
$$ where, in the last term, we have used identity (\ref{bi1}) and
eliminated $A_{dbc}$ via [C]. The last term in this expression
obviously vanishes identically. On the other hand the first term also
vanishes, because of identity (\ref{bi3}). Thus a simple consequence
of equation [C] is $C_{dabc}\nabla^{a}K^{d}=0.$ Thus, since the metric
is generic (in fact for this result we only need that it is
$\Lambda^2$-generic), we can conclude that
$$
\nabla^{[a}K^{d]}=0.
$$
Therefore, at least locally, there
exists a function $\Up$ such
that 
\begin{equation}
K_d=\nabla_d\Up. \label{grad}
\end{equation}
Thus, we have shown that our conditions [C], [B]
are equivalent to the necessary conditions (\ref{1ic}), (\ref{2ic}) for
a metric to be conformally Einstein.\\

To prove the sufficiency we first take $\nabla^c$ of [C]. This, after
using the identity (\ref{bi1}) and the definition of the Bach tensor
(\ref{bach}), takes the form 
$$B_{ab}+\Rho^{dc}C_{dabc}-C_{dabc}\nabla^cK^d+(n-3)K^d K^c C_{dabc}=0.$$ 
Now, subtracting from this equation our second condition [B] we get 
\begin{equation}\label{c3barer}
C_{dabc}(\Rho^{dc}-\nabla^cK^d+K^dK^c)=0.
\end{equation}
Next we differentiate equation [C] and skew to obtain 
$$
\nabla_{a_1}A_{ca_2a_3}-C_{a_2a_3cd}\nd_{a_1}K^d -K^d\nd_{a_1}C_{a_2a_3cd}=0.
$$
Then using \nn{ABianchiid}, the Weyl Bianchi identity \nn{WBianchiid}, 
and [C] once more we obtain
$$
C_{a_2a_3cd}(\Rho_{a_1}{}^d -\nd_{a_1}K^d+K_{a_1}K^d)=0
$$
or equivalently
\begin{equation}\label{c3*barer}
C^*_{b_1\cdots b_{n-2}cd}(\Rho^{b_1}{}^d -\nd^{b_1}K^d+K^{b_1}K^d)=0
\end{equation}
But this condition and \nn{c3barer} together imply that
$\Rho^{dc}-\nabla^cK^d+K^dK^c$ must be a pure trace, due to \nn{c3bar}
and \nn{c3*bar}. Thus,
$$\Rho^{dc}-\nabla^cK^d+K^dK^c=\frac{1}{n}Tg^{cd}.$$
This, when compared with our previous result (\ref{grad}) on $K^a$, and
with the conformal Einstein equations (\ref{confein}), shows that our
metric can be scaled to the Einstein metric with the function $\Up$
defined by (\ref{grad}). This proves the following theorem.
\begin{theorem}\label{I}
A generic metric $g_{ab}$ on an $n$-manifold $M$ 
is conformally Einstein if and only if its
Cotton tensor $A_{abc}$ and its Bach tensor $B_{ab}$ satisfy 
$$
\begin{array}{ll}
A_{abc}+K^d C_{dabc}=0 & \quad {\rm [C]}\\
&\\
B_{ab}+(n-4)K^dK^cC_{dabc}=0 & \quad {\rm [B]}
\end{array}
$$
for some vector field $K^a$ on $M$.
\end{theorem}
We will show below, and in the next section that [C] is
conformally invariant and that, while [B] is not conformally
invariant, the system [C], [B] is.  In particular [B] is conformally
invariant for metrics satisfying [C], the conformal C-space metrics.
Next note that, although we settled dimension 3 earlier, the above
theorem also holds in that case since the Weyl tensor vanishes
identically and the Bach tensor is just a divergence of the Cotton
tensor.  In other dimensions we can easily eliminate the {\it
undetermined} vector field $K^d$ from this theorem. Indeed, using the
tensor $||C||^{-1}\tilde{C}_{ed}^{~~bc}$ of (\ref{til}) and applying it on the
condition [C] we obtain
$$||C||^{-1}\tilde{C}_{ed}^{~~bc}A_{abc}+\frac{1}{2}(K_e g_{da}-K_dg_{ea})=0.$$
By contracting over the indices $\{ea\}$, this gives 
\begin{equation}\label{Kform}
K^d=\frac{2}{1-n}||C||^{-1}\tilde{C}^{dabc}A_{abc}.
\end{equation}

Inserting \nn{Kform} into the equations [C] and [B] of  
Theorem \ref{I}, 
we may reformulate the theorem as 
 the observation that 
a generic metric $g_{ab}$ on a $n$-manifold $M$ (where $n\geq 4$)  
is conformally Einstein if and only if its
Cotton tensor $A_{abc}$ and its 
Bach tensor $B_{ab}$ satisfy 
$$~~~~~~~~~(1-n)A_{abc}+2||C||^{-1}C_{dabc}\tilde{C}^{defg}A_{efg}=0~~~~~~~~~~~~~~~~~~~~~[C']
$$
and 
$$(n-1)^2B_{ab}+4(n-4)
||C||^{-2}\tilde{C}^{defg}C_{dabc}\tilde{C}^{chkl}A_{efg}A_{hkl}=0.~~~~~~~~[B']
$$
These are equivalent to conditions polynomial in the curvature. 
Multiplying the left
hand sides of [$\mbox{C}'$] and [$\mbox{B}'$] by, 
respectively, $||C||$ and $||C||^2$ we
obtain natural (pseudo-)Riemannian invariants which are obstructions
to a metric being conformally Einstein,
$$
F^{1}_{abc}:= (1-n)||C||A_{abc}+2C_{dabc}\tilde{C}^{defg}A_{efg}
$$
and 
$$
F^{2}_{ab}= (n-1)^2||C||^2 B_{ab}+4(n-4)
\tilde{C}^{defg}C_{dabc}\tilde{C}^{chkl}A_{efg}A_{hkl} .
$$ 
By construction the first of these is conformally covariant (see
below), the second tensor is conformally covariant for metrics such
that $F^{1}_{abc}=0$, and we have the following theorem.
\begin{theorem}\label{II}
A generic metric $g_{ab}$ on an $n$-manifold $M$ (where $n\geq 4$)  
is conformally Einstein if and only if the natural invariants $F^{1}_{abc}$ and 
$F^{2}_{ab}$ both vanish.
\end{theorem}

\noindent
{\bf Remarks}
\begin{itemize}
\item In dimension $n=4$ there exist examples of metrics satisfying
  the Bach equations [B] and not being conformally Einstein (see
  e.g. \cite{PN}). In higher dimensions it is straightforward to write
  down generic Riemannian metrics which, at least at a formal level, have
  vanishing Bach tensor but for which the Cotton tensor is
  non-vanishing. Thus the integrability condition [B] does
  not suffice to guarantee the conformally Einstein property of the
  metric. In Section \ref{examples} we discuss an example of special
  Robinson-Trautman metrics, which satisfy the condition [C] and do
  not satisfy [B]. (These are generic.) Thus condition [C] alone is not sufficient to
  guarantee the conformal Einstein property.
\item The development above parallels and generalises the tensor
treatment in \cite{CTP} which is based in dimension 4.  It should be
pointed out however that there are some simplifications in dimension
4.  Firstly $ F^{2}_{ab}$ simplifies to $9||C||^2 B_{ab}$. It is thus
sensible to use the conformally invariant Bach tensor $B_{ab}$ as a
replacement for $F^2$ in dimension 4. Also note, from the
development in \cite{CTP}, that the conditions that a metric $g_{ab}$
be generic may be characterised in a particularly simple way in
Lorentzian dimension four.  In this case they are equivalent to the
non-vanishing of at least one of the following two quantities,
$$C^3:=C_{abcd}C^{cd}_{~~ef}C^{efab}~~~~{\rm or}~~~~  
*C^3:=*C_{abcd}*C^{cd}_{~~ef}*C^{efab},$$ where
$*C_{abcd}=C^*_{abcd}=\vol_{abef}C^{ef}_{~~cd}$.
\end{itemize}

\subsection{Conformal invariants giving a sharp obstruction}\label{newinvts}

We will show in the next section that the system [C], [B] has a
natural and valuable geometric interpretation. However its value, or
the equivalent obstructions $F^1$ and $F^2$, as a test for conformally
Einstein metrics is limited by the requirement that the metric is
generic. Many metrics fail to be generic. For example in the setting
of dimension 4 Riemannian structures any selfdual metric fails to be
generic (and even fails to be $\Lambda^2$-generic), since any
anti-selfdual two form is a solution of \nn{c3}; at each point the
solution space of \nn{c3} is at least three dimensional (see section
4.3 for an explicit Ricci-flat example of this type).  In the
remainder of this section we show that there are natural conformal
invariants that are more effective, for detecting conformally Einstein
metrics, than the pair $F^1$ and $F^2$.

Let us say that a (pseudo-)Riemannian manifold is {\em weakly generic}
if the only solution $V^d$ to
\begin{equation}\label{weak}
C_{abcd}V^d=0
\end{equation}
is $V^d=0$.  From \nn{til} it is immediate that all $\Lambda^2$-generic spaces are
weakly generic and hence all generic spaces are weakly generic.
Via elementary multilinear algebra one can show that 
on weakly generic manifolds there is a tensor
field $\tilde{D}^{ab}{}_c{}^d$ with the property that
$$
\tilde{D}^{ac}{}_d{}^e C_{bc}{}^d{}_e=-\delta^a_b.
$$
Of course $\tilde{D}^{ab}{}_c{}^d$ is not uniquely determined by this
property. However in many settings there is a canonical choice. For
example in the case of Riemannian signature $g$ is weakly generic if
and only if $L^a_b:= C^{acde}C_{bcde}$ is invertible. Let us write
$\tilde{L}^a_b$ for the tensor field which is pointwise adjugate of
$L^a_b$.  $\tilde{L}^a_b$ is given by a formula which is
a partial contraction polynomial (and homogeneous of degree $2n-2$) in
the Weyl curvature and for any structure we have
$$
\tilde{L}^a_b L^b_c =||L||\delta^a_c ,
$$ where $||L||$ denotes the determinant of $L^a_b$.
 Let us define 
$$
D^{acde}:=- \tilde{L}^a_b C^{bcde} 
$$ Then $D^{acde}$ is a natural conformal covariant defined on all
structures. On weakly generic Riemannian structures, or
pseudo-Riemannian structures where we have $||L||\neq 0$,  there is a 
canonical choice for $\tilde{D}$, viz.\
\begin{equation}\label{Dtil}
\tilde{D}^{acde}:=||L||^{-1}D^{acde}= - ||L||^{-1}\tilde{L}^a_b C^{bcde}.
\end{equation}
On the other hand if  $g$ is $\Lambda^2$-generic
we may take 
\begin{equation}\label{Lam2D}
\tilde{D}^{acde}:= \frac{2}{1-n}||C||^{-1}\tilde{C}^{acde}.
\end{equation}
as was done implicitly in the previous section. Recall
$\tilde{C}^{ac}{}_d{}^e$ is conformally invariant and natural. The
examples \nn{Dtil} and \nn{Lam2D} are particularly important since
they are easily described and apply to any dimension (greater than
3). However in a given dimension there are many other possibilities
which lead to formulae of lower polynomial order if we know, or are
prepared to insist that, certain invariants are non-vanishing (see
\cite{Edgar} for a discussion in the context of $\Lambda^2$-generic
structures). For example in the setting of dimension 4 and Lorentzian
signature, $\Lambda^2$-generic implies
$C^3=C_{ab}{}^{cd}C_{cd}{}^{ef}C_{ef}{}^{ab}\neq 0$ and then one may
take $\tilde{D}^{acde}= C^{de}{}_{fg}C^{fgca}/C^3$ cf.\ \cite{CTP}.
In any case let us fix some choice for $\tilde{D}$.  Note that since
the Weyl curvature $C_{bc}{}^d{}_e$ for a metric $g$ is the same as
the Weyl tensor for a conformally related metric $\wh g$, it follows
that we can (and will) use the same tensor field
$\tilde{D}^{ab}{}_c{}^d $ for all metrics in the conformal class.

For weakly generic manifolds it is straightforward to give a
conformally invariant tensor that vanishes if and only if the manifold
is conformally Einstein. For the remainder of this section we 
assume the manifold is weakly generic.

We have observed already that the conformally Einstein manifolds are a
subclass of conformal C-spaces. Recall that a conformal C-space is a
(pseudo-)Riemannian manifold which admits a 1-form field $K_a$ which
solves the equation [C]:
$$
A_{abc}+K^d C_{dabc}=0 .
$$
If $K_1^d$ and $K_2^d$ are both solutions to [C] then, evidently,
$(K_1^d- K_2^d) C_{dabc}=0$. Thus, if the manifold is weakly generic,
 $K_1^d=K_2^d$.  In fact if $K_d$ is a solution to [C] then
clearly
\begin{equation}\label{formK}
K_d=\tilde{D}_d{}^{abc}A_{abc} ,
\end{equation}
which also shows that at most one vector field $K^d$ solves [C] on
weakly generic manifolds.  From either result, combined with the
observations that the Cotton tensor is preserved by constant conformal
metric rescalings and that constant conformal rescalings take Einstein
metrics to Einstein metrics, gives the following results.
\begin{proposition}\label{weakunique}
On a manifold with a weakly generic metric $g$, the equation [C] has
at most one solution for the vector field $K^d$.  

Either there are no metrics, conformally related
to $g$, that have vanishing Cotton tensor or the space of such metrics
is one dimensional.
Either there are no Einstein metrics, conformally
related to $g$, or the space of such metrics is one dimensional. 
\end{proposition}
\noindent If $g$ is a metric with vanishing Cotton tensor we will say 
this is a {\em C-space scale}.

Now, for an alternative view of conformal C-spaces, we may take \nn{formK} as the
{\em definition} of $K_d$. Note then that from \nn{Rhotrans}, a
routine calculation shows that $\wh{A}_{abc}=A_{abc}+\Up^kC_{kabc}$,
and so (using the conformal invariance of $\tilde{D}_d{}^{abc}$)
$K_d=\tilde{D}_d{}^{abc}A_{abc}$ has the conformal transformation
$$
\widehat{K}_d=K_d-\Up_d,
$$ where $\wh{A}_{abc}$ and $\wh{K}_d$ are calculated in terms of the
metric $\wh g= e^{2\Up}g$ and $\Up_a=\nd_a \Up$.  Thus $A_{abc}+K^d
C_{dabc}$ is conformally invariant.  From
proposition \ref{weakunique} and \nn{formK} this tensor is a {\em
sharp obstruction} to conformal C-spaces in the following sense.
\begin{proposition}\label{Cob}
A weakly generic manifold is a conformal C-space if 
and only if the conformal invariant 
$$
 A_{abc}+\tilde{D}^{dijk}A_{ijk}C_{dabc}
$$
vanishes.
\end{proposition}
In any case where $\tilde{D}^{dijk}$ is given by a Riemannian
invariant formulae rational in the curvature and its covariant
derivatives (e.g. $g$ is of Riemannian signature, or that $g$ is
$\Lambda^2$-generic) we can multiply the invariant here by an
appropriate polynomial invariant to obtain a natural conformal
invariant. Indeed, in the setting of $\Lambda^2$-generic metrics, the
invariant $F^1_{abc}$ (from section \ref{KNTstyle}) is an
example. Since, on $\Lambda^2$-generic manifolds, the vanishing of
$F^1_{abc}$ implies that \nn{Kform} is locally a gradient, we have the
following theorem.
\begin{theorem}\label{L2Cotton}
For a $\Lambda^2$-generic Riemannian or pseudo-Riemannian metric $g$ the conformal
covariant $F^1_{abc}$,
$$
(1-n)||C||A_{abc}+2C_{dabc}\tilde{C}^{defg}A_{efg}
$$
vanishes if and only if $g$ is conformally related to a \underline{Cotton metric} (i.e. a metric $\wh{g}$ such that its Cotton tensor vanishes,  $\wh{A}_{abc}=0$). 
\end{theorem}
\noindent In the case of Riemannian signature $\Lambda^2$-generic metrics we may
replace the conformal invariant $F^1_{abc}$ in the theorem with the 
conformal invariant, 
\begin{equation}\label{RL2Cotton}
||L||A_{abc}- C^{efgh}A_{fgh}\tilde{L}^d_e C_{dabc} \quad n\geq 4.
\end{equation}
In dimension 4 there is an even simpler invariant. 
Note that in dimension 4  we have 
\begin{equation}\label{4id}
4C^{abcd}C_{abce}=|C|^2\delta^d_e
\end{equation}
 where $|C|^2:=C^{abcd}C_{abcd} $ and so $L$ is a multiple of the
 identity. Eliminating, from \nn{RL2Cotton}, the factor of
 $(|C|^2)^3$ and a numerical scale we obtain the conformal invariant
$$
|C|^2A_{abc}- 4C^{defg}A_{efg}C_{dabc} \quad n=4,
$$ which again can be used to replace $F^1_{abc}$ in the theorem for dimension 4
$\Lambda^2$-generic metrics.

We can also characterise conformally Einstein spaces.
\begin{proposition}\label{Eprop}
 A weakly generic metric $g$ is conformally Einstein if and only if
 the conformally invariant tensor
$$ E_{ab}:=\mbox{\rm Trace-free}\big[\Rho_{ab}-\nd_a(
\tilde{D}_{bcde}A^{cde})+\tilde{D}_{aijk}A^{ijk} \tilde{D}_{bcde}A^{cde} \big]
$$
vanishes.
\end{proposition}
\noindent{\bf Proof:}
The proof that $E_{ab}$ is conformally invariant is a simple
calculation using \nn{Rhotrans} and the transformation formula for
$K_d =\tilde{D}_d{}^{abc}A_{abc}$.

If $g$ is conformally Einstein then there is a gradient $\Upsilon_a$ such that 
$$
\mbox{Trace-free}\big[\Rho_{ab}-\nd_a\Up_b+ \Up_a\Up_b \big] =0.
$$
From Section \ref{necessary} this implies $\Up_a$ solves the C-space equation 
(see \nn{1ic}) and hence, from \nn{formK}, $\Up_a = \tilde{D}_{aijk}A^{ijk}$, and so $E_{ab}=0$.

Conversely suppose that $E_{ab}=0$. Then the skew part of $E_{ab}$
vanishes and since $\Rho_{ab}$ and $\tilde{D}_{aijk}A^{ijk} \tilde{D}_{bcde}A^{cde}$
are symmetric we conclude that $\tilde{D}_{bcde}A^{cde}$ is closed and hence,
locally at least, is a gradient. 
\quad $\Box$

Now suppose  $||L||\neq 0$ and take $\tilde{D}_{abcd}$ to be given as in  \nn{Dtil}.
Note that since $E_{ab}$ is conformally invariant it follows that  $||L||^2 E_{ab}$ is
conformally invariant.
This expands to 
$$
G_{ab}:= \mbox{Trace-free}\big[||L||^2 \Rho_{ab}-||L||\nd_a(
D_{bcde}A^{cde})+\quad\quad\quad\quad\quad\quad\quad\quad\quad\quad\quad\quad $$
$$\quad\quad\quad\quad\quad\quad\quad\quad\quad\quad\quad\quad(\nd_a||L||)(
D_{bcde}A^{cde}) +D_{aijk}A^{ijk} D_{bcde}A^{cde} \big].
$$ This is natural by construction. Since it is given by a universal
polynomial formula which is conformally covariant on structures for
which $||L||$ is non-vanishing, it follows from an elementary
polynomial continuation argument that it is conformally covariant on
any structure. Note $||L||$ is a conformal covariant of weight
$-4n$. Thus we have the following theorem on manifolds of dimension $n\geq 4$.
\begin{theorem}\label{nob}
The natural invariant $G_{ab}$ is a conformal covariant of weight $-8n$.
  A manifold with a weakly generic Riemannian metric $g$ is conformally Einstein if
  and only if $G_{ab}$ vanishes. The same is true on pseudo-Riemannian manifolds where 
the conformal invariant $||L||$ is non-vanishing.
\end{theorem}

Recall that in dimension 4 we have the identity \nn{4id}. Thus
 $||L||\neq 0$ if and only if $|C|^2\neq0$ and we obtain a
 considerable simplification. In particular the invariant $G_{ab}$ has
 an overall factor of $(|C|^2)^6$ that we may divide out and still
 have a natural conformal invariant.  This corresponds to taking
 $(|C|^2)^2E_{ab}$ with $\tilde{D}^{abcd}=-\frac{4}{|C|^2}C^{abcd}$.
 Hence we have a simplified obstruction as follows.
\begin{theorem}\label{4ob}
The natural invariant 
$$ \mbox{\rm Trace-free}\big[\quad\quad\quad\quad\quad\quad\quad\quad\quad\quad\quad\quad\quad\quad\quad\quad\quad\quad\quad\quad\quad\quad\quad\quad\quad\quad\quad\quad\quad\quad$$
$$  (|C|^2)^2\Rho_{ab}+ 4
|C|^2 \nd_a (C_{bcde}A^{cde}) - 4 C_{bcde}A^{cde} \nd_a |C|^2 + 
16C_{aijk}A^{ijk} C_{bcde}A^{cde}\big]$$ 
is conformally covariant of weight $-8$. 

A 4-manifold with $|C|^2$ nowhere vanishing
is conformally Einstein if and only if this invariant vanishes.
\end{theorem}
 In the case of Riemannian 4-manifolds, requiring $|C|^2$ non-vanishing
 is the same as requiring the manifold to be weakly generic. In
 this setting this is a very mild assumption; note that $|C|^2=0$ at
 $p\in M$ if and only if $C_{abcd}=0$ at $p$ (and so the manifold is
 conformally flat at $p$).

 Note also that
if we denote by $F_{ab}$ the natural invariant in the Theorem then on
Riemannian 4 manifolds the (conformally covariant) scalar function
$F_{ab}F^{ab}$ is an equivalent sharp obstruction to the manifold
being conformally Einstein.

Now suppose we are in the setting of $\Lambda^2$-generic
structures (of any fixed signature). Then $E_{ab}$ is well defined and
conformally invariant with $\tilde{D}_{abcd}$ given by
\nn{Lam2D}. Thus again by polynomial continuation we can conclude that
the natural invariant obtained by expanding  $||C||^2 E_{ab}$, viz.\
$$
\bar{G}_{ab}:=\quad\quad\quad\quad\quad\quad\quad\quad\quad\quad\quad\quad\quad\quad\quad\quad\quad\quad\quad\quad\quad\quad\quad\quad\quad\quad\quad\quad\quad\quad\quad\quad\quad\quad$$
$$ \mbox{Trace-free}\big[ (1-n)^2||C||^2 \Rho_{ab}-2(1-n)||C||\nd_a(
\tilde{C}_{bcde}A^{cde})+\quad\quad\quad\quad\quad\quad\quad$$
$$\quad\quad\quad\quad 2(1-n)(\nd_a||C||)(
\tilde{C}_{bcde}A^{cde}) +4\tilde{C}_{aijk}A^{ijk} \tilde{C}_{bcde}A^{cde} \big].
$$ 
 is conformally covariant on any
structure (i.e.\ not necessarily $\Lambda^2$-generic).
Thus we have the following theorem on manifolds of dimension $n\geq 4$.
\begin{theorem}\label{Lam2ob}
The natural invariant $\bar{G}_{ab}$ is a conformal covariant of weight $2n(1-n)$.
  A manifold with a $\Lambda^2$-generic metric $g$ is conformally Einstein if
  and only if $\bar{G}_{ab}$ vanishes. 
\end{theorem}

We should point out that there is further scope, in each specific
dimension, to obtain simplifications and improvements to Theorems
\ref{nob} and \ref{Lam2ob} along the lines of Theorem \ref{4ob}. For
example in dimension 4 the complete 
contraction $C^3=C_{ab}{}^{cd}C_{cd}{}^{ef}C_{ef}{}^{ab}$, 
mentioned earlier, is a
conformal covariant which is independent of $|C|^2$ (see e.g.\
\cite{Pen8}). Thus on pseudo-Riemannian structures this may be
non-vanishing when $|C|^2=0$. There is the identity
$$
4 C_{jb}{}^{cd}C_{cd}{}^{ef}C_{ef}{}^{ib}= \delta^i_jC_{ab}{}^{cd}C_{cd}{}^{ef}C_{ef}{}^{ab}
$$ and this may be used to construct a formula for $\tilde{D}$ (and
then $K_d$ via \nn{Kform}) alternative to \nn{Dtil}
and \nn{Lam2D}. (See \cite{CTP} for this and some other examples.) 

Finally note that although generally we need to make some restriction
on the class of metrics to obtain a canonical formula for
$\tilde{D}_{bcde}$ in terms of the curvature, in other circumstances it
is generally easy to make a choice and give a description of a
$\tilde{D}$. For example in a non-Riemannian setting one can calculate
in a fixed local basis field and artificially nominate a Riemannian
signature metric. Using this to contract indices of the Weyl curvature
(given in the set basis field) one can then use the formula for $L$
and then $D$. In this way Proposition \ref{Eprop} is an effective and practical
means of testing for conformally Einstein metrics, among the class
weakly generic metrics, even when it does not lead to a natural invariant. 

\section{A geometric derivation and new obstructions}\label{tractor}

The derivation of the system of theorem \ref{I} appears ad hoc.  We
will show that in fact [C] and [B] are two parts (or components) of a
single conformal equation that has a simple and clear geometric
interpretation. This construction then easily yields new obstructions. 
This is based on the observation that conformally
Einstein manifolds may be characterised as those admitting a parallel
section of a certain vector bundle. The vector bundle concerned is the
(standard) conformal tractor bundle. This bundle and its canonical
conformally invariant connection are associated structures for the
normal conformal Cartan connection of \cite{C}. The initial development
of the calculus associated to this bundle dates back to the work of
T.Y. Thomas \cite{T} and was reformulated and further developed in a
modern setting in \cite{BEGo}. For a comprehensive treatment exposing
the connection to the Cartan bundle and relating the conformal case to
the wider setting of parabolic structures see
\cite{Cap-Gover,Cap-Gover2}.  The calculational techniques,
conventions and notation used here follow \cite{GP} and \cite{Goadv}.

\subsection{Conformal geometry and tractor calculus}\label{tractorsect}

We first introduce some of the basic objects of conformal tractor
calculus.  It is useful here to make a slight change of point of
view. Rather than take as our basic geometric structure a Riemannian
or pseudo-Riemannian structure we will take as our basic geometry only
a conformal structure. This simplifies the formulae involved and their
conformal transformations. It is also a conceptually sound move since
conformally invariant operators, tensors and functions are exactly the
(pseudo-)Riemannian objects that descend to be well defined objects on
a conformal manifold. A signature $(p,q)$ {\em conformal structure}
$[g]$ on a manifold $M$, of dimension $n\ge 3$, is an equivalence
class of metrics where $\wh{g}\sim g$ if $\wh g =e^{2\Upsilon} g$ for
some $\Upsilon\in \ce$.  A conformal structure is equivalent to a ray
subbundle $\cq$ of $S^2T^*M$; points of $\cq$ are pairs $(g_x,x)$
where $x\in M$ and $g_x$ is a metric at $x$, each section of $\cq$
gives a metric $g$ on $M$ and the metrics from different sections
agree up to multiplication by a positive function.  The bundle $\cq$
is a principal bundle with group $\bR_+$, and we denote by $\ce[w]$
the vector bundle induced from the representation of $\bR_+$ on $\bR$
given by $t\mapsto t^{-w/2}$.  Sections of $\ce[w]$ are called a {\em
conformal densities of weight $w$} and may be identified with
functions on $\cq$ that are homogeneous of degree $w$, i.e., $f(s^2
g_x,x)=s^w f(g_x,x)$ for any $s\in \bR_+$.  We will often use the same
notation $\ce[w]$ for the space of sections of the bundle. Note that
for each choice of a metric $g$ (i.e., section of $\cq$, which we term
a {\em choice of conformal scale}), we may identify a section $f\in
\ce[w]$ with a function $f_g$ on $M$ by $f_g(x)=f(g_x,x)$. This
function is conformally covariant of weight $w$ in the sense of
Section \ref{tensor}, since if $\widehat{g}=e^{2\Up}g$, for some $\Up
\in \ce$, then $f_{\widehat{g}}(x)=f(e^{2\Up_x}
g_x,x)=e^{w\Up_x}f(g_x,x)=e^{w\Up_x}f_g(x)$.  Conversely conformally
covariant functions determine homogeneous sections of $\cq$ and so
densities.  In particular, $\ce[0]$ is canonically identified with
$\ce$.

Note that there is a tautological function $\bg$ on $\cq$ taking
values in $S^2T^*M$.  It is the function which assigns to the point
$(g_x,x)\in \cq$ the metric $g_x$ at $x$.  This is homogeneous of
degree 2 since $\bg (s^2 g_x,x) =s^2 g_x$. If $\xi$ is any positive
function on $\cq$ homogeneous of degree $-2$ then $\xi \bg$ is
independent of the action of $\bR_+$ on the fibres of $\cq$, and so
$\xi \bg$ descends to give a metric from the conformal class. Thus
$\bg$ determines and is equivalent to a canonical section of
$\ce_{ab}[2]$ (called the conformal metric) that we also denote $\bg$
(or $\bg_{ab}$). This in turn determines a canonical section
$\bg^{ab}$ (or $\bg^{-1}$) of $\ce^{ab}[-2]$ with the property that
$\bg_{ab}\bg^{bc}=\delta_a^c$ (where $\delta_a{}^c $ is kronecker
delta, i.e., the section of $\ce^c_a$ corresponding to the identity
endomorphism of the the tangent bundle).  In this section the
conformal metric (and its inverse $\bg^{ab}$) will be used to raise
and lower indices. This enables us to work with density valued
objects. Conformally covariant tensors as in section \ref{tensor}
correspond one-one with conformally invariant density valued tensors.
Each non-vanishing section $\si$ of $\ce[1]$ determines a
metric $g^\si$ from the conformal class by
\begin{equation}\label{cscale}
g^\si:=\si^{-2}\bg. 
\end{equation}
Conversely if $g\in [g]$ then there is an up-to-sign unique
$\si\in\ce[1]$ which solves $g=\si^{-2}\bg$, and so $\si$ is termed a
choice of conformal scale.  Given a choice of conformal scale, we
write $ \nabla_a$ for the corresponding Levi-Civita connection. For
each choice of metric there is also a canonical connection on $\ce[w]$
determined by the identification of $\ce[w]$ with $\ce$, as described
above, and the exterior derivative on functions. We will also call
this the Levi-Civita connection and thus for tensors with weight,
e.g.\ $v_{a}\in\ce_{a}[w]$, there is a connection given by the Leibniz
rule. With these conventions the Laplacian $ \Delta$ is given by
$\Delta=\bg^{ab}\nd_a\nd_b= \nd^b\nd_b\,$.

We next define the standard tractor bundle over $(M,[g])$.
It is a vector bundle of rank $n+2$ defined, for each $g\in[g]$,
by  $[\ce^A]_g=\ce[1]\oplus\ce_a[1]\oplus\ce[-1]$. 
If $\wh g=e^{2\Up}g$, we identify  
 $(\alpha,\mu_a,\tau)\in[\ce^A]_g$ with
$(\wh\alpha,\wh\mu_a,\wh\tau)\in[\ce^A]_{\wh g}$
by the transformation
\begin{equation}\label{transf-tractor}
 \begin{pmatrix}
 \wh\alpha\\ \wh\mu_a\\ \wh\tau
 \end{pmatrix}=
 \begin{pmatrix}
 1 & 0& 0\\
 \Up_a&\delta_a{}^b&0\\
- \tfrac{1}{2}\Up_c\Up^c &-\Up^b& 1
 \end{pmatrix} 
 \begin{pmatrix}
 \alpha\\ \mu_b\\ \tau
 \end{pmatrix} .
\end{equation}
It is straightforward to verify that these identifications are
consistent upon changing to a third metric from the conformal class,
and so taking the quotient by this equivalence relation defines the
{\em standard tractor bundle} $\ce^A$ over the conformal manifold.
(Alternatively the standard tractor bundle may be constructed as a
canonical quotient of a certain 2-jet bundle or as an associated
bundle to the normal conformal Cartan bundle \cite{Cap-Gover2}.)  The
bundle $\ce^A$ admits an invariant metric $ h_{AB}$ of signature
$(p+1,q+1)$ and an invariant connection, which we shall also denote by
$\nabla_a$, preserving $h_{AB}$.  In a conformal scale $g$, these are
given by
$$
 h_{AB}=\begin{pmatrix}
 0 & 0& 1\\
 0&\bg_{ab}&0\\
1 & 0 & 0
 \end{pmatrix}
\text{ and }
\nabla_a\begin{pmatrix}
 \alpha\\ \mu_b\\ \tau
 \end{pmatrix}
 =
\begin{pmatrix}
 \nabla_a \alpha-\mu_a \\
 \nabla_a \mu_b+ \bg_{ab} \tau +\Rho_{ab}\alpha \\
 \nabla_a \tau - \Rho_{ab}\mu^b  \end{pmatrix}. 
$$
It is readily verified that both of these are conformally well-defined,
i.e., independent of the choice of a metric $g\in [g]$.  Note that
$h_{AB}$ defines a section of $\ce_{AB}=\ce_A\otimes\ce_B$, where
$\ce_A$ is the dual bundle of $\ce^A$. Hence we may use $h_{AB}$ and
its inverse $h^{AB}$ to raise or lower indices of $\ce_A$, $\ce^A$ and
their tensor products.

In computations, it is often useful to introduce 
the `projectors' from $\ce^A$ to
the components $\ce[1]$, $\ce_a[1]$ and $\ce[-1]$ which are determined
by a choice of scale.
They are respectively denoted by $X_A\in\ce_A[1]$, 
$Z_{Aa}\in\ce_{Aa}[1]$ and $Y_A\in\ce_A[-1]$, where
 $\ce_{Aa}[w]=\ce_A\otimes\ce_a\otimes\ce[w]$, etc.
 Using the metrics $h_{AB}$ and $\bg_{ab}$ to raise indices,
we define $X^A, Z^{Aa}, Y^A$. Then we
immediately see that 
$$
Y_AX^A=1,\ \ Z_{Ab}Z^A{}_c=\bg_{bc}
$$
and that all other quadratic combinations that contract the tractor
index vanish. This is summarised in Figure~\ref{XYZfigure}. 
\begin{figure}
$$
\begin{array}{l|ccc}
& Y^A & Z^{Ac} & X^{A}
\\
\hline
Y_{A} & 0 & 0 & 1
\\
Z_{Ab} & 0 & \delta_{b}{}^{c} & 0
\\
X_{A} & 1 & 0 & 0
\end{array}
$$
\caption{Tractor inner product}
\label{XYZfigure}
\end{figure}

It is clear from \eqref{transf-tractor} that the first component 
$\alpha$ is independent of the choice of a representative $g$ and 
hence $X^A$ is conformally invariant. 
For $Z^{Aa}$ and $Y^A$, we have the transformation laws:
  \begin{equation}\label{XYZtrans}
  \wh Z^{Aa}=Z^{Aa}+\Up^aX^A, \quad
  \wh Y^A=Y^A-\Up_aZ^{Aa}-\frac12\Up_a\Up^aX^A.
\end{equation}

Given a  choice
of conformal  scale we have the corresponding Levi-Civita connection
on tensor and density bundles.  In this setting we can use the coupled 
Levi-Civita tractor connection to act on sections of the tensor product 
of a tensor bundle with a tractor bundle. This is defined by the Leibniz 
rule in the usual way.  For example if
$ u^b V^C \alpha\in \ce^b\otimes \ce^C\otimes \ce[w]=:\ce^{bC}[w]$
then 
$\nd_a u^b V^C \alpha = (\nd_a u^b) V^C \alpha +
u^b(\nd_a V^C) \alpha + u^b V^C \nd_a \alpha$.
Here $\nd$ means the Levi-Civita
connection on $ u^b\in \ce^b$ and $ \alpha\in \ce[w]$,
while it denotes the tractor
connection on $ V^C\in \ce^C$. In particular with this convention we have 
\begin{equation}\label{connids}
\nd_aX_A=Z_{Aa},\quad
\nd_aZ_{Ab}=-\Rho_{ab}X_A-Y_A\bg_{ab}, 
\quad \nd_aY_A=\Rho_{ab}Z_A{}^b.
\end{equation}

Note that if $V$ is a section of $ \ce_{A_1\cdots A_\ell}[w]$,
then the coupled Levi-Civita tractor
connection on $V$ is not conformally invariant but transforms just as the
Levi-Civita connection transforms on densities of the same weight:
$\wh{\nd}_a V = \nd_a V + w\Up_a V$.

Given a choice of conformal scale, the {\em tractor-$D$ operator} 
$$
D_A\colon\ce_{B \cdots E}[w]\to\ce_{AB\cdots E}[w-1]
$$
is defined by 
\begin{equation}\label{Dform}
D_A V:=(n+2w-2)w Y_A V+ (n+2w-2)Z_{Aa}\nabla^a V -X_A\Box V, 
\end{equation} 
 where $\Box V :=\Delta V+w \Rho_b{}^b V$.  This also turns out to be
 conformally invariant as can be checked directly using the formulae
 above (or alternatively there are conformally invariant constructions
 of $D$, see e.g.\ \cite{Gosrni}).

The curvature $ \Omega$ of the tractor connection 
is defined by 
\begin{equation}\label{curvature}
[\nd_a,\nd_b] V^C= \Omega_{ab}{}^C{}_EV^E 
\end{equation}
for $ V^C\in\ce^C$.  Using
\eqref{connids} and the usual formulae for the curvature of the
Levi-Civita connection we calculate (cf. \cite{BEGo})
\begin{equation}\label{tractcurv}
\Omega_{abCE}= Z_C{}^cZ_E{}^e C_{abce}-2X_{[C}Z_{E]}{}^e A_{eab}
\end{equation}

From the tractor curvature we obtain a related higher order conformally 
invariant curvature quantity by the formula (cf.\ \cite{Gosrni,Goadv})
$$
W_{BC}{}^E{}_F:=
\frac{3}{n-2}D^AX_{[A} \Omega_{BC]}{}^E{}_F .
$$
It is straightforward to  verify that this can be re-expressed as
follows,
\begin{equation}
\label{Wform}
W_{ABCE}=(n-4)Z_A{}^aZ_B{}^b\Omega_{abCE}
-2X_{[A}Z_{B}{}_{]}^b\nd^p \Omega_{pbCE}.
\end{equation}
This tractor field has an important relationship to the ambient metric
of Fefferman and Graham.  For a conformal manifold of signature
$(p,q)$ the ambient manifold \cite{FeffGr} is a signature $(p+1,q+1)$
pseudo-Riemannian manifold with $\cq$ as an embedded
submanifold. Suitably homogeneous tensor fields on the ambient
manifold upon restriction to $\cq$ determine tractor fields on the
underlying conformal manifold \cite{CapGoamb}.  In particular, in
dimensions other than 4, $W_{ABCD}$ is the tractor field equivalent
to $(n-4){\aR}|_\cq$ where $\aR$ is the curvature of
the Fefferman-Graham ambient metric.

\subsection{Conformally Einstein manifolds}\label{ceinm}

Recall that we say a Riemannian or pseudo-Riemannian metric $g$ is
conformally Einstein if there is a scale $\Up$ such that the Ricci
tensor, or equivalently the Schouten tensor, is pure trace. Thus we
say that a conformal structure $[g]$ is conformally Einstein if there
is a metric $\widehat{g}$ in the conformal class (i.e. $\widehat{g}\in
[g]$) such that the Schouten tensor for $\widehat{g}$ is pure trace.
We show here that a conformal manifold $(M,[g])$ is
conformally Einstein if and only if it admits a parallel standard
tractor $\Pa^A$ which also satisfy the condition that $X_A\Pa^A$ is nowhere
vanishing. Note that in a sense the ``main condition'' is that $\Pa$ is parallel since  
 $X_A\Pa^A\neq 0$ is an open condition.
 In more detail we have the following result.
\begin{theorem}\label{cein}
On a conformal manifold $(M,[g])$ there is a 1-1 correspondence
between conformal scales $\si\in \ce[1]$, such that
$g^\si=\si^{-2}\bg$ is Einstein, and parallel standard tractors $\Pa$
with the property that $X_A\Pa^A$ is nowhere vanishing. The mapping
from Einstein scales to parallel tractors is given by 
$\si\mapsto \frac{1}{n}D_A \si$ while the inverse is $\Pa^A \mapsto X_A\Pa^A$.
\end{theorem}
\noindent{\bf Proof:}
Suppose that $(M,[g])$ admits a parallel standard tractor $\Pa^A$ such
that $\si:=X_A\Pa^A$ is nowhere vanishing. Since $\si\in \ce[1]$ and
is non-vanishing it is a conformal scale. Let $g$ be the metric from
the conformal class determined by $\si$, that is $g=g^\si=\si^{-2}\bg$ as in \nn{cscale}.
In terms of the tractor bundle splitting determined by this metric $\Pa^A $ is given by
some triple with $\si$ as the leading entry, $[\Pa^A]_{g}=(\si,\mu_a,\tau)$. 
From the formula for the invariant
connection we have \renewcommand{\arraystretch}{1}
\begin{equation}\label{ndP}
0=[\nabla_a \Pa^B]_g =
\left(\begin{array}{c} \nabla_a \si-\mu_a \\
                       \nabla_a \mu_b+ \bg_{ab} \tau +\V_{ab}\si \\
                       \nabla_a \tau - \V_{ab}\mu^b  \end{array}\right) . 
\end{equation}
\renewcommand{\arraystretch}{1.5} Thus $\mu_a= \nd_a \si$, but $ \nd_a
\si=0$ by the definition of $\nd$ in  the scale $\si$. 
Thus $\mu_a$ vanishes, and the second tensor  equation from \nn{ndP} simplifies to
$$
\V_{ab}\si =- \bg_{ab} \tau ,
$$ 
showing that the metric $g$ is Einstein. Note that tracing the
display gives $\tau=-\frac{1}{n}\J \si$.

To prove the converse let us now suppose that $\si$ is a conformal scale so that
$g=\si^{-2}\bg$ is an Einstein metric. That is,
 for this metric $g$,
$\V_{ab}$ is pure trace. Let us work in this conformal scale. Then we have
$\V_{ab} =\frac{1}{n}\bg_{ab}\J $. Thus $\nd^a\V_{ab}=(1/n) \nd_b \J$.
On the other hand comparing this to the contracted Bianchi identity
 $\nd^a\V_{ab}= \nd_b \J$ 
we have that $\nd_a \J =0$.
Now, we define a tractor field $\Pa^A$ by $\Pa^A:= \frac{1}{n}D^A \si$.
Then $[\Pa]_{g^\si}:=(\si, 0 , -\frac{1}{n}\J \si)$. Consider the tractor 
connection on this.
We have
$$
[\nd_a \Pa^B]_{g}=\left(\begin{array}{c} \nabla_a \si \\
                        -\frac{1}{n}\bg_{ab} \J\si  +\V_{ab}\si \\
                       -\frac{1}{n}(\si \nabla_a \J+ \J \nd_a \si)   \end{array}\right) . 
$$
Once again, by the definition of the Levi-Civita connection $\nd$ as
determined by the scale $\si$, we have $\nd \si=0$. Since $\V_{ab}
=\frac{1}{n}\bg_{ab}\J $ the second entry also vanishes.  The last
component also vanishes from $\nd \J=0$ and $\nd \si=0$. So $\Pa$ is a
parallel standard tractor satisfying $X_A\Pa^A=\si\neq 0$.
\quad $\Box$\\
\noindent{\bf Remarks:} 
\begin{itemize}
\item Note that $h(\Pa,\Pa)$ is a conformal
invariant of density weight 0. In fact from the formulae above, in the
Einstein scale, $h(\Pa,\Pa)=-\frac{2}{n}\si^2 \J$. Recall that in this section
$\J=\bg^{ab}\Rho_{ab}$ and so has density weight $-2$ and 
$$
\si^2 \J=\si^2\bg^{ab}\Rho_{ab} =g^{ab}\Rho_{ab}.
$$ 
That is $-\frac{n}{2}h(\Pa,\Pa)$ is the trace of Schouten tensor using the
metric determined by $\si$.  Since $\nd$ preserves the tractor metric
and $\Pa$ is parallel we recover the (well known) result that
$\Rho_{ab}$ (and its trace) is constant for Einstein metrics.

\item Suppose we drop the condition $\si:= X^A\Pa^a\neq 0$. If $\Pa^A$
  is parallel then from \nn{ndP} it follows that $\mu_a=\nd_a \si $.
  Furthermore tracing the middle entry on the right-hand-side of
  \nn{ndP} implies that $\tau= -\frac{1}{n}\Box \si$. Thus if $\nd_a
  \Pa_B=0$ at $p\in M$ then at $p$ we have $\Pa_B =\frac{1}{n}D_B
  \si$. So for parallel $\Pa^A$, $X_A \Pa^A$ vanishes on a
  neighbourhood if and only if $\Pa^A$ vanishes on the same
  neighbourhood.  If $\Pa^A$ is non-zero then $X_A\Pa^A$ may vanish on
  submanifolds of dimension at most $n-1$. The points of these
  submanifolds are conformal singularities  for the metric $g=\si^{-2}\bg$.

\item On dimension 4 spin manifolds it is straightforward to show that
the standard tractor bundle is isomorphic to the second exterior
power of Penrose's \cite{ot} local twistor bundle.  Under this
isomorphism $\Pa$ may be identified with the {\em infinity twistor}
(defined for spacetimes). The relationship to conformal Einstein manifolds
is well known \cite{MB,FS} in that setting. 

\item We should also point out that the theorem above can
alternatively be deduced, via some elementary arguments but without any calculation, 
from the construction
of the tractor connection as in \cite{BEGo}. 
\end{itemize}

Next we make some elementary observations concerning parallel tractors.
\begin{lemma}\label{DWN}
On a conformal manifold let $N$ be a parallel section of the standard
tractor bundle $\bT$. Then: \newline
$$
\Omega_{bc}{}^D{}_EN^E =0 \quad \mbox{ and } \quad W_{BCDE}N^E =0 
$$
\end{lemma}
\noindent{\bf Proof:}
By assumption we have  $\nd_a N^D=0$. Thus $\Omega_{bc}{}^D{}_EN^E=[\nd_b,\nd_c]N^D=0$ and 
the first result is established.

Next $W_{A_1A_2}{}^D{}_E N^E= \frac{3}{n-2}(D^{A_0}X_{A_0}
Z_{A_1}{}^bZ_{A_2}{}^c\Omega_{bc}{}^D{}_E)N^E$, where, as usual,
sequentially labelled indices e.g.\ $A_0,A_1,A_2$ are implicitly skewed
over. Now the quantity $X_{A_0}
Z_{A_1}{}^bZ_{A_2}{}^c\Omega_{bc}{}^D{}_E$ has (density) weight $-1$,  so
from the formula  \nn{Dform} for $D$, we have
$$ 
\begin{array}{lll}
\lefteqn{(D^{A_0}X_{A_0}
Z_{A_1}{}^bZ_{A_2}{}^c\Omega_{bc}{}^D{}_E)N^E =}&&
\\
&&
(4-n)Y ^{A_0}X_{A_0}
Z_{A_1}{}^bZ_{A_2}{}^c\Omega_{bc}{}^D{}_E N^E
\\
&&
+(n-4)(Z^{A_0a}{}\nd_aX_{A_0}
Z_{A_1}{}^bZ_{A_2}{}^c\Omega_{bc}{}^D{}_E)N^E
\\
&&
- (X^{A_0}\Delta X_{A_0}
Z_{A_1}{}^bZ_{A_2}{}^c\Omega_{bc}{}^D{}_E) N^E
\\
&&
+ \J X^{A_0} X_{A_0}
Z_{A_1}{}^bZ_{A_2}{}^c\Omega_{bc}{}^D{}_E N^E ,
\end{array}
$$ where $\nd$ and $\Delta$ act on everything to their right within
the parentheses.  The first and last terms on the right-hand-side
vanish from the previous result.  (In fact for last term we could also
use that $X^{A_0}X_{A_0} Z_{A_1}{}^bZ_{A_2}{}^c=0$.) Next observe
that, since $\nd N=0$, we have $$(Z^{A_0a}{}\nd_aX_{A_0}
Z_{A_1}{}^bZ_{A_2}{}^c\Omega_{bc}{}^D{}_E)N^E =
Z^{A_0a}{}\nd_a(X_{A_0}
Z_{A_1}{}^bZ_{A_2}{}^c\Omega_{bc}{}^D{}_EN^E)=0$$ where we have again
used the earlier result, $\Omega_{bc}{}^D{}_EN^E =0 $.  Similarly
$$(X^{A_0}\Delta X_{A_0} Z_{A_1}{}^bZ_{A_2}{}^c\Omega_{bc}{}^D{}_E)
N^E= X^{A_0}\Delta (X_{A_0} Z_{A_1}{}^bZ_{A_2}{}^c\Omega_{bc}{}^D{}_E
N^E ) =0.  $$ \quad $\Box$

From the Lemma it follows immediately that on conformally Einstein
manifolds the parallel tractor $\Pa$, of Theorem \ref{cein}, satisfies
$\Omega_{bc}{}^D{}_E\Pa^E =0$ and $ W_{BCDE}\Pa^E =0$.  In general the
converse is also true. More accurately we have the result given in the
following theorem.  Before we state that, note that since the Weyl
curvature is conformally invariant it follows that the equations
\nn{c3}, \nn{c3bar} and \nn{c3*bar} are conformally invariant. Thus if
any metric from a conformal class is generic then all metrics from the
class are generic and we will describe the conformal class as generic.
\begin{theorem} \label{maintr}
A generic conformal manifold of dimension $n\neq 4$ is conformally Einstein if and only if
there exists a  tractor field $\Pa^A\in \ce^A$ such that $X_A\Pa^A\neq 0$ and 
$$
W_{BCDE}\Pa^E =0.
$$
A generic conformal manifold of dimension $n=4 $ is conformally Einstein if and only if
there exists a  tractor field $\Pa^A\in \ce^A$ such that $X_A\Pa^A\neq 0$,
$$
\Omega_{bc}{}^D{}_E\Pa^E =0 \quad \mbox{and} \quad W_{BCDE}\Pa^E =0.
$$
\end{theorem}
\noindent{\bf Proof:} We have shown that on a conformally Einstein
manifold there is a (parallel) standard tractor field satisfying\\
 \IT{i}  $X_A\Pa^A \neq 0$,\\ 
\IT{ii} $\Omega_{bc}{}^D{}_E\Pa^E =0$, \\
  \IT{iii} $ W_{BCDE}\Pa^E =0$ .\\

It remains to prove the relevant converse statements.
First we observe that given  \IT{i}, \IT{ii} is exactly the
conformal C-space equation. 
From above we have that 
$$ 
\Omega_{abCE}= Z_C{}^cZ_E{}^e C_{abce}-X_{C}Z_{E}{}^e A_{eab}+ X_{E}Z_{C}{}^e A_{eab} 
$$
A general tractor $\Pa^A\in \ce^A$  may
be expanded to  
$$
\Pa^E =  Y^E\si + Z^{Ed}\mu_d+X^E\t,
$$
where $\si=X_A\Pa^A$ and we assume this is non-vanishing. 
Hence 
\begin{equation}\label{OmPa}
\Omega_{abCE}\Pa^E= \si Z_{C}{}^c A_{cab}  +Z_C{}^c \mu^d C_{abcd}- X_{C}\mu^d A_{dab} .
\end{equation}
Setting this to zero, as required by \IT{ii}, implies that the coefficient of 
$Z_C{}^c$ must vanish, i.e., 
$
 \si A_{cab} +\mu^d C_{abcd} =0 ,
$
or 
\begin{equation}\label{Cspace}
A_{cab} +K^d C_{dcab} =0, \quad K^d:=-\si^{-1}\mu^d ,
\end{equation}
 which is exactly the conformal C-space equation [C] as in theorem
\ref{I}.  
Contracting this with $\mu^c$ (or $K^c$) annihilates the second
term and so
$$
\mu^d A_{dab}=0,
$$
whence the coefficient of $X_C$ in \nn{OmPa} vanishes as a
consequence of the earlier equation and it is shown that (with \IT{i})
$\Omega_{abCE}\Pa^E=0$ is exactly the conformal C-space equation.

Now recall
$$
W_{BCDE}= (n-4)Z_B{}^bZ_C{}^c\Om_{bcDE}-2 X_{[B}Z_{C]}{}^c \nd^a\Om_{acDE},
$$
and so, in dimensions other 4, $W_{BCDE}\Pa^E=0$ implies
$\Om_{bcDE}\Pa^E=0$ (and hence the conformal C-space equation). From
the display we see that $W_{BCDE}\Pa^E=0$ also implies that $
\Pa^E\nd^a\Om_{acDE}=0$ or equivalently $
\si^{-1}\Pa^E\nd^a\Om_{acDE}=0$. Once again using the formulae for the
tractor connection we obtain
\begin{equation}\label{divtrc}
\nd^a\Om_{acDE}= 
(n-4)Z_D{}^dZ_E{}^e A_{cde}- X_{D}Z_{E}{}^e B_{ec} + X_{E}Z_{D}{}^e B_{ec}
\end{equation}
where $B_{ec}$
is the Bach tensor. 
Hence  
$\si^{-1}\Pa^E\nd^a\Om_{acDE}=0$ expands to
$$
-(n-4)Z_D{}^d K^e A_{cde}+ X_{D} K^e B_{ec}
+ 
Z_{D}{}^d B_{dc} =0 .
$$
From the coefficient of $Z_{D}{}^d $ we have 
$$
B_{dc}-(n-4) K^e A_{cde} =0
$$
which, with the conformal C-space equation (and since $B$ is
symmetric), gives
\begin{equation}\label{B}
B_{cd} + (n-4) K^e K^a C_{acde} =0
\end{equation}
which is exactly the second equation [B] of Theorem 1. If this holds
then it follows at once that $K^c B_{cb}=0$ and so in 
the expansion of  $\si^{-1}\Pa^E\nd^a\Om_{acDE}=0$
 the coefficient of $X_D$ vanishes without further restriction.
Thus we have shown that in dimensions other than 4 the single
conformally invariant tractor equation $W_{BCDE}\Pa^E=0$ is equivalent
to the two equations [C] and [B]. In dimension 4 it is clear from
\nn{Wform} that $W_{BCDE}\Pa^E=0$ is equivalent to $\Pa^E
\nd^a\Om_{acDE}=0$ and this with $\Pa^E \Om_{acDE}=0$ gives the pair
of equations [B],[C].  In either case then the theorem here now
follows immediately from Theorem \ref{I}.  \quad $\Box$\\
\noindent{\bf Remarks:} \begin{itemize}

\item Note that conditions \IT{i}, \IT{ii}
and \IT{iii}, as in the theorem, do not imply that $\Pa$ is parallel.
On the other hand the theorem shows that if there exists a standard
tractor $\Pa$ satisfying these conditions then (on generic manifolds)
also there exists a parallel standard tractor $\Pa'$ satisfying these
conditions.  Calculating in an Einstein scale, it follows from the
conformal C-space equation that one has $Z_A{}^a\Pa^A=Z_A{}^a\Pa^A=0$.
Hence that $\Pa'=f\Pa+\rho X$ for some section $\rho$ of $\ce[-1]$ and
non-vanishing function $f$. 
\item   Recall that in Section \ref{tractorsect} we pointed out that 
in dimensions other than 4,
$W_{ABCD}$ is the tractor field equivalent \cite{CapGoamb} to
$(n-4){\aR}|_\cq$ where $\aR$ is the curvature of the Fefferman-Graham
ambient metric. Thus, in these dimensions, the condition
$W_{ABCD}\Pa^D=0 $ is equivalent to the existence of a suitably
homogeneous and generic ambient tangent vector field along $\cq$ in
the ambient manifold which annihilates the ambient curvature.
\item  We had already observed in
section \ref{newinvts} that $A_{abc}+K^d C_{dabc}$ is conformally
invariant if we assume that $K_d$ has the conformal transformation law
$\wh{K}_a=K_a-\Up_a$ (where $\wh{g}=e^{2\Up}g$). From the proof above
we see this transformation formula fits naturally into the tractor picture and arises from 
\nn{transf-tractor} since $K_a$ is a density multiple of the middle
component of a tractor field according to \nn{Cspace}.  
\end{itemize}

\subsection{Sharp obstructions via tractors}

Theorem \ref{maintr} gives a simple interpretation of
Theorem \ref{I} in terms of tractor bundles. In the proof of this above, this
connection was made by recovering the familiar tensor equations from
section \ref{tensor}.  Here we first observe that entire derivation of
Theorem \ref{I} and its proof reduces to a few key lines if we work in
the tractor picture. This then leads to a stronger theorem as below.

 We
summarise the background first. From Theorem \nn{cein} we know that
the existence of a conformal Einstein structure is equivalent to the
existence of a parallel tractor $\Pa$ (at points where $X_A\Pa^A\neq
0$). This immediately implies that the tractor curvature $\Om_{abCD}$
satisfies
$$
\begin{array}{cl}
\Pa^D\Om_{abCD}=0 & \mbox{[$\tilde{\rm C}$]} \\
\Pa^D\nd^a\Om_{abCD}=0 & \mbox{[$\tilde{\rm B}$]}. 
\end{array}
$$ We have labelled these [$\tilde{\rm C}$] and [$\tilde{\rm B}$]
since (as shown in the proof above) the first equation is equivalent
to the earlier [C] and, given this, the second equation is equivalent
to the earlier equation [B]. The conformal invariance of the system
$[C],[B]$ is now immediate in all dimensions from the observation that
the conformal transformation of $\nd^a\Om_{abCD}$ is 
\begin{equation}\label{divtrctrans}
\wh{\nd^a\Om_{abCD}} = \nd^a\Om_{abCD} +(n-4) \Up^a\Om_{abCD},
\end{equation}
and whence
 the conformal transformation of the left-hand-side of equation \mbox{[$\tilde{\rm B}$]} is 
 $$
 \wh{\Pa^D\nd^a\Om_{abCD}} = \Pa^D\nd^a\Om_{abCD} +(n-4)
 \Up^a\Pa^D\Om_{abCD},
$$
where $\wh{g}=e^{2\Up}g$; from this it is immediate that
\mbox{[$\tilde{\rm B}$]} is invariant on metrics that
solve \mbox{[$\tilde{\rm C}$]}. We should point out that in dimension 4
it follows immediately from \nn{divtrc} that $\Pa^D\nd^a\Om_{abCD}=0
\Leftrightarrow \nd^a\Om_{abCD}=0 \Leftrightarrow B_{ab}=0$.  

Now we are interested in the converse. We will show that if the
displayed equations \mbox{[$\tilde{\rm C}$]} and \mbox{[$\tilde{\rm
B}$]} hold for some tractor $\Pa$ satisfying $X_A\Pa^A\neq 0$ then the
structure is conformally Einstein. Here is an alternative proof of
Theorem \ref{maintr} (and hence an alternative proof of Theorem
\ref{I}).  Equation [$\tilde{\rm C}$] implies that $\nd_{a_1}(
\Om_{a_2a_3CD}\Pa^D)=0$, where as usual sequentially labelled indices
are skewed over. From the Bianchi identity for the tractor curvature,
$\nd_{a_1} \Om_{a_2a_3CD}=0$, it follows that
\begin{equation}\label{*ver}
\Om_{a_2a_3CD}\nd_{a_1}\Pa^D=0 .
\end{equation}
Now equation [$\tilde{\rm C}$] implies [C], viz.\ $A_{cab} +K^d
C_{dcab} =0$.  As we saw earlier this (using that the metric is
$\Lambda^2$-generic) implies that $K_a$ is a gradient and that there
is a conformal scale such that the Cotton tensor $A_{cab}$
vanishes. In this special C-space scale (see Section \ref{newinvts})
it is clear that $K_a$ is also zero and \nn{*ver} simplifies (using
\nn{ndP} and \nn{tractcurv}) to $ \Rho_{a_1}{}^dC_{a_2a_3 cd}Z_{C}{}^c
=0 $ or equivalently
\begin{equation}\label{two}
C^*_{b_1\cdots b_{n-2}cd}\Rho^{b_1d}=0 .
\end{equation}
Note that if $C^*$ is suitably generic this already implies that
the metric that gives the special C-space scale is Einstein. 

Using only the weaker assumption that the manifold is generic in the
sense of section \ref{KNTstyle} we must also use [$\tilde{\rm B}$].
The argument is similar to the above. 
Equation [$\tilde{\rm C}$] implies 
$\nd^a(\Pa^D\Om_{abCD})=0$. 
Thus using  [$\tilde{\rm B}$] we have 
$$
(\nd^a \Pa^D)\Om_{abCD}=0 .
$$
In the special C-space scale this expands to 
$
\Rho^{ad}C_{ab cd}Z_{C}{}^c =0,
$
which is equivalent to 
\begin{equation}\label{one}
\Rho^{ad}C_{ab cd}=0. 
\end{equation}
Clearly equations \nn{one} and \nn{two} imply that $\Rho$ is pure
trace on generic manifolds and so the theorem is proved. In fact these equations \nn{one} and \nn{two} 
are respectively equations \nn{c3barer} and \nn{c3*barer} both 
written in the C-space scale. 

\vspace{3mm}

The construction of the system [$\tilde{\rm B}$] and [$\tilde{\rm C}$]
immediately suggests alternative systems. In particular we have the
following results which only requires the manifold to be 
weakly generic.
\begin{theorem} \label{tracchar} A weakly generic conformal manifold is conformally 
Einstein if and only if there exists a non-vanishing tractor field $\Pa^A\in \ce^A$ 
such that
$$
\begin{array}{cl}
\Pa^E\Om_{bcDE}=0 & \mbox{\rm [$\tilde{\rm C}$]} \\
\Pa^E\nd_a\Om_{bcDE}=0 & \mbox{\rm [$\tilde{\rm D}$]}. 
\end{array}
$$
The system \mbox{\rm [$\tilde{\rm C}$]}, \mbox{\rm [$\tilde{\rm D}$]}
 is conformally invariant. 
\end{theorem}
\noindent{\bf Proof:} Note that from \nn{conntrans}, and the
invariance of the tractor connection, we have
 $$
 \wh{\Pa^E\nd_a\Om_{bcDE}} = \Pa^E\nd_a\Om_{bcDE} -2
 \Up_a\Pa^E\Om_{bcDE} -\Up_b\Pa^E\Om_{acDE} 
- \Up_c\Pa^E\Om_{baDE} \quad \quad \quad \quad$$
$$\hspace*{6cm} +
\bg_{ab}\Up^k\Pa^E\Om_{kcDE}+ \bg_{ac} \Up^k\Pa^E\Om_{bkDE},
$$ 
where $\wh{g}=e^{2\Up}g$, and so [$\tilde{\rm D}$] is conformally
invariant if the conformally invariant equation [$\tilde{\rm C}$] is
satisfied; the system [$\tilde{\rm C}$], [$\tilde{\rm D}$] is 
conformally invariant.
 
If the manifold is conformally Einstein then there is a parallel
tractor $\Pa^E$. We have observed earlier that this satisfies
[$\tilde{\rm C}$]. Differentiating [$\tilde{\rm C}$] and then using
once again that $\Pa^E$ is parallel shows that [$\tilde{\rm D}$] is
satisfied.

Now we assume that [$\tilde{\rm C}$] and [$\tilde{\rm D}$] hold. If $
\Pa^E = Y^E\si + Z^{Ed}\mu_d+X^E\t, $ then $\Omega_{abCE}\Pa^E$ is
given by \nn{OmPa}.  Suppose that $X_A\Pa^A=\si=0$.  Then from
\nn{OmPa} we have $ \mu^d C_{abcd}=0$ (and $\mu^dA_{dab}=0$) and so, since the conformal
class is weakly generic, $\mu^d=0$. Thus $\Pa^E=\t X^E$ and
[$\tilde{\rm D}$] becomes $X^E\nd_a\Om_{bcDE}=0$. But, $\nd_a X^E=Z^E{}_a$ and from
\nn{tractcurv} $X^E\Om_{bcDE}=0$,  and so $
Z_D{}^d C_{bcda}-X_{D} A_{abc}= Z^E{}_a\Om_{bcDE}=0$. But this means
$C_{bcda}=0$ which contradicts the assumption that the conformal class
is weakly generic. So $X_A\Pa^A\neq 0$.

Now, differentiating  [$\tilde{\rm C}$] and then using [$\tilde{\rm D}$] we obtain
$$
\Om_{bcDE}\nd_a\Pa^E=0 .
$$ But, since the manifold is weakly generic, $\Om_{bcDE}$ must have
rank at least $n$ as a map $\Om_{bcDE}:\ce^{bcD}\to \ce_E$. Also, from
\nn{tractcurv} and [$\tilde{\rm C}$], $X^E$ and $\Pa^E$ are orthogonal
to the range. So the display implies that
$$
\nd_a\Pa^E= \alpha_a \Pa^E +\beta_a X^E ,
$$
for some 1-forms $\alpha_a$ and $\beta_a$.
 (An alternative explanation is to note, as earlier, that if $U^E$
is not a multiple of $X^E$ and $\Om_{bcDE} U^E=0$ then from
\nn{tractcurv} it follows that $U^E$ determines a non-trivial solution
of the equation [C]. Since $\Pa^E$ also determines such a solution it follows at once from 
Proposition \ref{weakunique}  that $U^E=\alpha \Pa +\beta X^E$.)
Differentiating again and alternating we obtain
$$
\Om_{ba}{}^E{}_D\Pa^D= 2\Pa^E\nd_{[b}\alpha_{a]} +2 \alpha_{[a}\alpha_{b]}\Pa^E
+2 \alpha_{[a}\beta_{b]} X^E + 2X^E\nd_{[b}\beta_{a]}+2\beta_{[a}Z^E{}_{b]} .
$$
The left-hand-side vanishes by assumption and of course $ \alpha_{[a}\alpha_{b]}\Pa^E=0 $. 
Contracting $X_E$ into the remaining terms brings us to 
$$
0=2\si \nd_{[a}\alpha_{b]} 
$$
and so $\alpha$ is closed. Locally then $\alpha_a=\nd_a f$ for some function f and so 
$\tilde{\Pa}^E:= e^{-f}\Pa^E$ satisfies
\begin{equation}\label{alpar}
\nd_a\tilde{\Pa}^E= \tilde{\beta}_a X^E
\end{equation}
for some 1-form $\tilde{\beta}_a$. Expanding $\tilde{\Pa}^E$: 
$\tilde{\Pa}^E= Y^E\tilde{\si} + Z^{Ed}\tilde{\mu}_d+X^E\tilde{\t}$ we have 
$0\neq X_E\tilde{\Pa}^e=\tilde{\si}$ and, from \nn{alpar}, the
equations 
$$
\begin{array}{l} \nabla_a \tilde{\si}-\tilde{\mu}_a =0 \\
                       \nabla_a \tilde{\mu}_b+ \bg_{ab} \tilde{\tau} +\V_{ab}\tilde{\si} =0 
   \end{array} 
$$ 
cf.\ \nn{ndP}. So for the metric $g:= \tilde{\si}^{-2}\bg$ we have
$\tilde{\mu}_a =\nabla_a \tilde{\si}=0 $ and $\V_{ab} + \bg_{ab}
\tilde{\tau}/\tilde{\si}=0$. That is the metric $g$ is Einstein (and
$\frac{1}{n}D_A \tilde{\si}$ is parallel).  \quad $\Box$

We have the following consequence of the theorem above.
\begin{corollary}\label{weakanysig}
A weakly generic pseudo-Riemannian or Riemannian metric $g$ on an $n$-manifold is 
conformally Einstein if and only if the natural invariants 
$$
\Om_{abK D_1}
 \cdots \Om_{cdLD_s} \nd_{e}\Omega_{fgPD_{s+1}} \cdots \nd_h \Omega_{k\ell QD_{n+2}} ,
$$ for $s=0,1,\cdots,n+1$, all vanish. Here the sequentially labelled
indices $D_1,\cdots , D_{n+2}$ are completely skewed over.
\end{corollary}
\noindent{\bf Proof:} 
The Theorem can clearly be rephrased to state that
 $g$ is conformally Einstein if and only if the map 
\begin{equation}\label{Ommap}
(\Om_{bcDE},\nd_a\Om_{bcDE}):\ce^{bcD}\oplus\ce^{abcD} \to \ce_E
\end{equation}
given by 
$$
(V^{bcD},W^{abcD})\mapsto V^{bcD} \Om_{bcDE}+W^{abcD}\nd_a\Om_{bcDE}
$$ fails to have maximal rank. But by elementary linear algebra this
happens if and only if the induced alternating multi-linear map to
$\Lambda^{n+2}(\ce^E)$ vanishes. This is equivalent to the claim in
the Corollary, since for any metric the tractor curvature satisfies
$\Om_{bcDE}X^E=0$.  \quad $\Box$ \\ 
If $M$ is oriented (which locally
we can assume with no loss of generality) then it is straightforward
to show that there is a canonical skew $(n+2)$-tractor consistent with
the tractor metric and the orientation. Let us denote this by
$\vol^{C_1 \cdots C_{n+2}}$. Using this, we could equally rephrase the
Corollary in terms of the invariants
$$
\vol^{D_1D_2\cdots D_sD_{s+1} \cdots D_{n+1}D_{n+2}}\Om_{abK D_1}
 \cdots \Om_{cdLD_s} \nd_{e}\Omega_{fgPD_{s+1}} \cdots \nd_h \Omega_{k\ell QD_{n+2}} ,
$$ for $s=0,1,\cdots,n+1$. These all vanish if and only if the metric is
conformally Einstein.

The natural invariants in the Lemma are given by mixed tensor-tractor
fields, rather pure tensors.  However by expanding $\Om_{abCD}$ and
$\nd_a\Om_{bcDE}$ using \nn{tractcurv} and \nn{connids} it is
straightforward to obtain an equivalent set of tensorial obstructions
from these. The system of obstructions so obtained is rather unwieldy
and could be awkward to apply in practise. Nevertheless this gives a
system of invariants, which works equally for all signatures.

As a final remark in this section we note that coming to Proposition
\ref{Eprop} via the tractor picture is also very easy. If we want to test whether
a scale $\si\in\ce[1]$ is an Einstein scale we define $\Pa_B:=
\frac{1}{n}D_B \si$ as in Theorem \ref{cein} and consider $\nd_a
\Pa_B$.  Calculating in terms of an arbitrary metric $g$ from the
conformal class we get $\nd_a\Pa_B= Z_B{}^b \si E_{ab}$, modulo terms
involving $X_B$, where $E_{ab}=\mbox{Trace-free}(\Rho_{ab}-\nd_a K_b
+K_aK_b )$ and $K_a:= -\si^{-1}\nd_a \si$.  Since $\si$ can only be an
Einstein scale if $\Omega_{bc}{}^D{}_E\Pa^E =0$ we obtain the
conformal C-space equation for $K_a$ and we are led to the conclusion
that the Riemannian invariant of the proposition is conformally
invariant and also the conclusion that it must vanish on conformal
Einstein manifolds.

\section{Examples}\label{examples}

Here we shed light on the various notions of generic metrics, mainly
by way of examples. First let us note that each of these is an open
condition on the moduli space of possible curvatures.  Thus in this
sense ``almost all'' metrics are generic (and hence
$\Lambda^2$-generic and weakly generic). The many components of the
Weyl curvature $C_{abcd}$ arise from a $\Lambda^2$-generic metric
unless they lie on the closed variety determined by the one 
condition $||C||=0$ where, recall, $||C||$ is the determinant of the
map \nn{Lam2map}. The metrics which fail to be weakly generic
correspond to a closed subspace contained in the $||C||=0$ variety.
In the Riemannian case this subvariety is given by $||L||=0$, where
recall $||L||$ is the determinant of $C^{acde}C_{bcde}$ and we show
below that in dimension 4 the containment is proper.

Another aim in this final section is to establish the independence of the conditions [C] and 
[B] from Section \ref{KNTstyle}. We assume that $n\geq 4$ throughout this section.

\subsection{Simple $n$-dimensional Robinson-Trautman metrics}
Let $Q$ be an $(n-2)$-dimensional space of constant curvature $\kappa$
and denote by $x^i$, $i=1,2,...,n-2$, standard stereographic coordinates on $Q$.   We
take $M={\bf R}^2\times Q$, with coordinates $(r,u,x^i)$, where
$(r,u)$ are coordinates along the ${\bf R}^2$, and equip $M$
with a subclass of Robinson-Trautman \cite{RobTra} metrics $g$ by 
\be
g=2\der u~[\der r+h(r)\der u~]~+~r^2~\frac{g_{ij}\der x^i\der
x^j}{(1+\frac{\kappa}{4}g_{kl}x^kx^l)^2}.\label{rt} 
\ee 
Here $g_{ij}={\rm diag}(\epsilon_1,\epsilon_2,...,\epsilon_{n-2})$,
$\epsilon_i=\pm 1$, $\kappa=1,0,-1$ and $h=h(r)$ is an arbitrary,
sufficiently smooth real function of variable $r$. In the following we
describe conformal properties of the metrics (\ref{rt}).\\

To calculate the Weyl tensor we introduce the null-orthonormal coframe
$(\theta^a)=(\theta^+,\theta^-,\theta^i)$ by
\be
\theta^+=\der u,~~~~\theta^-=\der r+h\der u,~~~~\theta^i=r\frac{\der x^i}{1+\frac{\kappa}{4}g_{kl}x^kx^l}.\label{ncof}
\ee 
In this coframe the metric takes the form $g=g_{ab}\theta^a\theta^b$ where
\be
g_{ab}=\begin{pmatrix} 
0&1&\\
1&0&\\
&&g_{ij}
\end{pmatrix}.
\ee
We lower and raise the indices by means of the matrix $g_{ab}$ 
and its inverse $g^{ab}$. The Levi-Civita connection 1-forms 
$$\G_{ab}=\G_{abc}\theta^c$$
are uniquely determined by \be \der
\theta^a+\G^a_{~b}\dz\theta^b=0~~~{\rm and}~~~\der
g_{ab}-\G_{ab}-\G_{ba}=0.  \ee Explicitly, we find that, the
connection 1-forms are \be
\G_{ij}=\frac{\kappa}{2r}(x_i\theta_j-x_j\theta_i),~~~~~~\G_{-j}=-\frac{1}{r}\theta_j,~~~~~~\G_{+j}=\frac{h}{r}\theta^j,~~~~~~\G_{+-}=h'\theta^+,
\ee where $h'=\frac{\der h}{\der r}$. (Observe that, due to the
constancy of the matrix elements of $g_{ab}$, the matrix $\G_{ab}$ is
skew, $\G_{ab}=-\G_{ba}$.)  The curvature 2-forms
$$
\Om_{ab}=\frac{1}{2}R_{abcd}\theta^c\dz\theta^d=\der\G_{ab}+\G_a^{~c}\dz\G_{cb}
$$
are
\be
\Om_{ij}=\frac{\kappa+2h}{r^2}\theta_i\dz\theta_j,~~~~~~\Om_{-j}=\frac{h'}{r}\theta^+\dz\theta_j,~~~~~~\Om_{+j}=\frac{h'}{r}\theta^-\dz\theta_j,~~~~~~\Om_{+-}=h''\theta^-\dz\theta^+,
\ee
with the remaining components determined by symmetry. The non-vanishing components of the Ricci tensor $$R_{ab}=R^c_{~acb}$$
and the Ricci scalar $$R=g^{ab}R_{ab}$$ are
\be
R_{ij}=[(n-3)\frac{\kappa+2h}{r^2}+\frac{2h'}{r}]g_{ij},~~~~~~~R_{+-}=(n-2)\frac{h'}{r}+h'',
\ee
$$R=(n-2)[(n-3)\frac{\kappa+2h}{r^2}+\frac{4h'}{r}]+2h''.$$
From this we conclude that metrics (\ref{rt}) are Einstein, 
$$R_{ab}=\Lambda g_{ab},$$ if and only if \be
h(r)=-\frac{\kappa}{2}+\frac{m}{r^{n-3}}+\frac{\Lambda}{2(n-1)}r^2,
\ee where $m$ and $\Lambda$ are constants. These metrics form the well
known $n$-dimensional Schwarzschild-(anti-)de Sitter 
2-parameter class in which $m$ is interpreted as the mass and
$\Lambda$ as the cosmological constant. 
(The space is termed de Sitter if $\Lambda>0$ and anti-de Sitter is $\Lambda<0$.)
Thus, we have the following
Proposition.  
\begin{proposition}~\\ The only Einstein metrics among the
Robinson-Trautman metrics
$$
g=2\der u~[\der r+h(r)\der u~]~+~r^2~\frac{g_{ij}\der x^i\der x^j}{(1+\frac{\kappa}{4}g_{kl}x^kx^l)^2}
$$
are the Schwarzschild-(anti-)de Sitter metrics, for which 
$$
h(r)=-\frac{\kappa}{2}+\frac{m}{r^{n-3}}+\frac{\Lambda}{2(n-1)}r^2.
$$ 
\end{proposition}
The Weyl tensor of metrics (\ref{rt}) has the
following non-vanishing components:
\be
C_{ijkl}=2\Psi (g_{ki}g_{jl}-g_{kj}g_{il}),~~~~~~C_{-i+k}=
(3-n)\Psi g_{ik},~~~~~~~C_{+-+-}=(3-n)(n-2)\Psi,
\ee
where 
$$
\Psi=\frac{1}{(n-1)(n-2)}~[~\frac{\kappa+2h}{r^2}-\frac{2h'}{r}+h''~],
$$
and the further non-vanishing components determined from these by the Weyl symmetries. 
Now, we consider the equation 
\be
C_{abcd}F^{cd}=0\label{gen1}
\ee
for the antisymmetric tensor $F_{ab}$. We easily find that 
$$C_{ijab}F^{ab}=4\Psi F_{ij},~~~~~C_{i+ab}F^{ab}=(3-n)\Psi
g_{ik}F^{k-},$$
$$C_{i-ab}F^{ab}=(3-n)\Psi g_{ik}F^{k+},~~~~~
C_{+-ab}F^{ab}=2(3-n)(n-2)\Psi F^{+-}.$$
Thus, if $\Psi\neq 0$, the equation (\ref{gen1}) has unique solution
$F_{ab}=0$. 
We pass to the equation 
\be
C_{abcd}H^{bd}=0\label{gen2}
\ee
for a symmetric and trace-free tensor $H_{ab}$. In the
null-orthonormal coframe (\ref{ncof}) the trace-free
condition reads
\be
H~+~2 H_{+-}=0,~~~~~~{\rm where}~~~~~H=g^{ik}H_{ik}.
\ee
Comparing this with
$$C_{ibkd}H^{bd}=2\Psi~[g_{ik}~(H+(3-n)H_{-+})~-~H_{ik}],$$
$$
C_{ib-d}H^{bd}=(n-3)\Psi g_{ik}H^{+k},~~~~~
C_{ib+d}H^{bd}=(n-3)\Psi g_{ik}H^{-k},$$
$$
C_{-b-d}H^{bd}=(n-2)(n-3)\Psi H^{++},~~~~~
C_{+b+d}H^{bd}=(n-2)(n-3)\Psi H^{--}$$
proves that the only solution of (\ref{gen2}) is $H_{ab}=0$. Thus we
have the following proposition.
\begin{proposition} \label{RTgeneric}
If
$$\Psi=\frac{1}{(n-1)(n-2)}~[~\frac{\kappa+2h}{r^2}-\frac{2h'}{r}+h''~]\neq
0$$ the Robinson-Trautman metrics 
$$g=2\der u~[\der r+h(r)\der u~]~+~r^2~\frac{g_{ij}\der x^i\der
  x^j}{(1+\frac{\kappa}{4}g_{kl}x^kx^l)^2}$$ are generic.
\end{proposition}
By a  straightforward calculation we obtain the following
proposition.
\begin{proposition}~\\
Each Robinson-Trautman metric 
for which $\Psi\neq 0$, 
satisfies the C-space condition [C] with a vector field $K_a$ given
by 
\be
K_a=\nabla_a~\log[~r^{\frac{(1-n)}{n-3}}\Psi^{\frac{1}{3-n}}~].\label{rozw}
\ee
\end{proposition}

From this and Propositions \ref{weakunique} and \ref{RTgeneric}
it follows that the
Robinson-Trautman metrics for which $\Psi\neq 0$ 
are conformal to Einstein metrics if and only if 
$$P_{ab}-\nabla_aK_b+K_a K_b-\frac{1}{n}(P-\nabla^c K_c+K^c
K_c)g_{ab}=0$$ with $K_a$ given by (\ref{rozw}). (Note that, by the
uniqueness asserted in Proposition \ref{weakunique}, this is equivalent 
to requiring  $E_{ab}=0$ with $E_{ab}$ as in Proposition \ref{Eprop}.) Inserting
$R_{ab}$ and $K_a$ into this equation one finds that the metric
(\ref{rt}) is conformal to an Einstein metric if and only if the
function $h=h(r)$ is given by
$$h(r)=-\frac{\kappa}{2}+\frac{m}{r^{n-3}}+\frac{\Lambda}{2(n-1)}
r^2.$$ 
This means that among the considered Robinson-Trautman metrics 
the only metrics which are conformal to Einstein metrics are those
belonging to the 2-parameter Schwarzschild-de Sitter
family. So we have  the
following conclusions. The Robinson-Trautman metrics (\ref{rt}):
\begin{itemize}
\item are all generic,
\item all satisfy C-space condition,
[C]
\item in general do not satisfy the Bach condition [B].
\end{itemize}
In fact from the conformal invariance of the system [C],[B] (see Section \ref{ceinm}) and
the condition of being generic, the same conclusions hold for all
metrics conformally related to Robinson-Trautman metrics.

This, when along with 4-dimensional examples of metrics satisfying the
Bach conditions [B] and not being conformal to Einstein \cite{BaiEast,PN}, proves
independence of the two conditions [C] and [B].

\subsection{$n$-dimensional pp-waves}

We noted in Section \ref{newinvts} that there are weakly generic
metrics that fail to be $\Lambda^2$-generic, and hence fail to be
generic. Here we observe that, there are conformally non-flat metrics
that fail to be even weakly generic.

We consider the $n$-dimensional pp-wave metric
$$g=2\der u~[\der r+h(x^i,u)\der u]~+~g_{ij}\der x^i\der x^j,$$
where $g_{ij}$ is are the components of a constant non-degenerate $(n-2)\times (n-2)$ matrix.
This, in the coframe 
$$
\theta^+=\der u,~~~~~~\theta^-=\der r+h\der u,~~~~~\theta^i=\der x^i ,
$$
has curvature forms
$$
\Om_{i+}=-h_{,ik}\theta^k\dz\theta^+,~~~~~~\Om_{ij}=\Om_{i-}=\Om_{+-}=0.
$$
So the Ricci scalar vanishes, $R=0$, and the only non-vanishing
components of the Ricci and the Weyl tensors are
$$R_{++}=-2g^{ij}h_{,ij},~~~~~~C_{i+j+}=\frac{2}{n-2}[g_{ij}~g^{kl}h_{,kl}-(n-2)h_{,ij}],$$
apart from the components determined by these via symmetries. 
Thus, this metric is Einstein if and only if the function $h=h(x^i,u)$ is harmonic in the $x^i$ variables, 
$$ g^{ij}h_{,ij}=0,$$
in which case it is also Ricci flat. Whether this is satisfied or not it is clear that 
 the vector field 
\be
K=f\partial_r,\label{przyl}
\ee 
where $f$ is any non-vanishing function, satisfies 
\be
C_{abcd}K^d=0.\label{rrr}
\ee 
Thus, the pp-wave metric is not weakly generic. It is worth noting that 
if the trace-free part of the matrix $h_{,ij}$ is invertible the vector 
(\ref{przyl}) is the most general solution of equation (\ref{rrr}). However, 
if it is not invertible, there are more vectors $K$ which satisfy (\ref{rrr}).

\subsection{4-dimensional hyperK\"ahler Ricci flat metrics}

Another interesting class of metrics that are weakly generic but not
$\Lambda^2$-generic or generic can be found in the complex setting.
Consider a 4-dimensional non-flat Riemannian metric, which is Ricci
flat and which admits three K\"ahler structures $I,J,K$ such that they
satisfy quaternionic identities, e.g. $IJ+JI=0$, $K=I J$. We claim
that all such manifolds are weakly generic, but not $\Lambda^2$-generic
\cite{PA}. To see this, first consider the Riemann tensor viewed as an
endomorphism $R(.): \La^2T^*M\to\La^2T^*M$.  Since the fundamental
forms $\om_I,\om_J,\om_K$, associated with $I,J,K$, are each parallel
we have $R(\om_I)=R(\om_J)=R(\om_K)=0$. On the other hand from
Ricci flatness we
have $R(.)=C(.)$, where $C(.)$ is the Weyl tensor,  also considered as and
endomorphism $C(.): \La^2T^*M\to\La^2T^*M $. 
Hence also
$C(\om_I)=C(\om_J)=C(\om_K)=0$, which means that the metric is not
$\Lambda^2$-generic. 

On the other hand if there existed a vector field $V$ such
that $C_{abcd}V^d=0$ then, because of the invariance property of $C$
with respect of the structures $I,J,K$ also $C_{abcd}(IV)^d$,
$C_{abcd}(JV)^d$ and $C_{abcd}(KV)^d$ would vanish.  Since on a hyperK\"ahler 
4-manifold a quadruple $(V,IV,JV,KV)$ associated with any  non-vanishing vector
$V$ constitutes a basis of vectors, at every point, we conclude that in
such a case $C_{abcd}$ (and therefore the Riemann tensor)
vanishes. Thus, if the manifold is not flat then we can
conclude that the Weyl tensor admits only $V=0$ as a
solution to $C_{abcd}V^d=0$.

Thus we have the following proposition
\begin{proposition}
Every non-flat 4-dimensional Riemannian Ricci flat hyperK\"{a}hler manifold
is weakly generic but not $\Lambda^2$-generic.
\end{proposition}

As a local example of this type we consider an open subset 
$$U=\{(z_1,z_2)\in{\bf C}^2~~{\rm such~that}~~z_1+\bar{z}_1-2z_2\bar{z}_2>0\}$$
of ${\bf C}^2$ equipped with the metric
$$g=\frac{1}{\sqrt{z_1+\bar{z}_1-2z_2\bar{z}_2}}(\der z_1-2\bar{z}_2\der z_2)
(\der \bar{z}_1-2z_2\der\bar{z}_2)+4\sqrt{z_1+\bar{z}_1-2z_2\bar{z}_2}\der z_2\der\bar{z}_2.$$
This metric is well known \cite{nurprzan} to be Ricci flat and hyperK\"{a}hler with the 
hyperK\"{a}hler structure given by
$$I=M\otimes n+\bar{M}\otimes \bar{n}-N\otimes m-\bar{N}\otimes\bar{m}$$
$$J=i(\bar{M}\otimes \bar{n}-M\otimes n+\bar{N}\otimes \bar{m}-N\otimes m)$$
$$K=i(\bar{M}\otimes\bar{m}-M\otimes m+N\otimes n-\bar{N}\otimes\bar{n}),$$
where
$$M=\frac{\der z_1-2\bar{z}_2\der z_2}{(z_1+\bar{z}_1-2z_2\bar{z}_2)^{\frac{1}{4}}},~~~~~
N=2(z_1+\bar{z}_1-2z_2\bar{z}_2)^{\frac{1}{4}}\der z_2,$$
and
$$m=(z_1+\bar{z}_1-2z_2\bar{z}_2)^{\frac{1}{4}}\partial_1,~~~~~
n=\frac{\partial_2+2\bar{z}_2\partial_1}{2(z_1+\bar{z}_1-2z_2\bar{z}_2)^{\frac{1}{4}}}.$$
Since
$$
|C|^2=\frac{24}{(z_1+\bar{z}_1-2z_2\bar{z}_2)^3}\neq 0,
$$
this metric is weakly generic. On the other hand,  from our considerations above, it is 
not $\Lambda^2$-generic and has its Weyl tensor is annihilated by the three fundamental 
K\"{a}hler forms 
$$\om_I=M\dz\bar{N}+\bar{M}\dz N$$
$$\om_J=i(\bar{M}\dz N-M\dz\bar{N})$$
$$\om_K=i(\bar{M}\dz M+N\dz\bar{N}).$$
It is worth noting that this metric also admits almost K\"ahler 
non-K\"ahler structure of opposite orientation than $I,J,K$ \cite{nurprzan}.

\end{document}